\documentclass[twoside]{article}
\usepackage{sfuEng}

\usepackage{graphicx}
\usepackage{color}

\usepackage{epsfig}
\usepackage{euscript}
\usepackage{textcomp}
\usepackage{amsfonts}
\usepackage{amssymb}
\usepackage{amstext}
\usepackage{amsthm}
\usepackage{newlfont}
\usepackage{subfigure}
\usepackage{epstopdf}

\tit{Dynamics of a symmetrically decoupled three-dimensional point transformation}

\abstract{In this work, we give some results obtained on the dynamics of a symmetrically decoupled three-dimensional point transformation. We are interested, in particular, in the study of its parametric plane and in its phase space by highlighting the existence of chaotic attractors.}


\begin{document}

\maketitBegin
\author{Hacene Gharout}{Laboratoire des Math\'ematiques Appliqu\'ees\\ Facult\'e des Sciences Exactes\\ Universit\'e de Bejaia\\  06000 Bejaia, ALGERIE}{gharouthacene@gmail.com}
\author{Abdel-kaddous Taha} {INSA, University of Toulouse\\ 135 Avenue de Rangueil \\31077 Toulouse, FRANCE}{taha@insa-toulouse.fr}
\author{Nourredine Akroune}{Laboratoire des Math\'ematiques Appliqu\'ees\\ Facult\'e des Sciences Exactes\\ Universit\'e de Bejaia\\  06000 Bejaia, ALGERIE}{akroune$\_$n@yahoo.fr}

\maketitEnd

\newcommand{\p}{\partial}
\renewcommand{\geq}{\geqslant}
\renewcommand{\leq}{\leqslant}
\renewcommand{\baselinestretch}{1}

\newtheorem{definition}{Definition}[section]
\newtheorem{remark}{Remark}
\newtheorem{propr}{Proprety}
\newtheorem{Proposition}{Proposition}
\newtheorem{theorem}{Theorem}[section]

\renewcommand{\thetheorem}{\arabic{section}.\arabic{theorem}}

\section{Introduction}
We will speak of a symmetrically decoupled three-dimensional point transformation, if this transformation can be put in the form:
\begin{equation}
  T \left\{
  \begin{array}{ll}
    x_{n+1} = f(y_{n}), \\
    y_{n+1} = g(z_{n}), \\
    z_{n+1} = h(x_{n}).
  \end{array}
\right.
\end{equation}
Such transformations are found in different fields of science, for example, in economics the duopoly of Cournot (1838) \cite{book1Cournot}; this problem has been the subject of a large number of publications T.Puu (2004) \cite{APuSU}, J.S.Canovas (2004) \cite{ACARU} and F.Tramontana et al in 2011 \cite{A29}. Another case that is modilized by such recurrences and that we meet in physics is that of systems with delay, see for example the work of L.Larger and D.Fournier-Prunaret (2005) \cite{ALArFOU} and L.Larger (2010) \cite{ALDU}. They are also found in the case of the security of the transmissions, see for example J.Xu, D.Fournier-Prunaret, A.K.Taha and P.Charg\'e (2010) \cite{A28} and H.Gharout et al. \cite{GHAROUT}. \\
We will focus here on the recurrence, defined by:
\begin{equation}\label{Eq301}
  T \left\{
  \begin{array}{ll}
    x_{n+1} = y_{n}, \\
    y_{n+1} = z_{n}, \\
    z_{n+1} = x^{2}_{n}+b.
  \end{array}
\right.
\end{equation}
\section{Properties of the map $T$}

We present in this section some definitions and properties of symmetrically decoupled point transformations, taken from A.Agliari et al. \cite{A27}, necessary to the understanding of this work. \\ The transformation $T$ is of the form:
\begin{equation}\label{Eq303}
  T \left\{
  \begin{array}{ll}
    x^{'} = f(y), \\
    y^{'} = g(z), \\
    z^{'} = h(x).
  \end{array}
\right.
\end{equation}
With $x=x_{n}$, $y=y_{n}$ and $z=z_{n}$; likewise $x^{'}=x_{n+1}$, $y^{'}=y_{n+1}$ and $z^{'}=z_{n+1}$.\\

The study of the cycles of the system defined by (\ref{Eq301}) will be done, by using the functions: $ H (x) = f (g (h (x)))$, $ F (y) = g (h (f (y))) $ and $ G (z) = h (f (g (z)))$. Using one-dimensional transformations $ H$, $ F $ and $ G $ it is possible to build the iterations of $ T $. Indeed:
\begin{eqnarray}
  T^{3k}(x,y,z)   &=& (H^{k}(x), F^{k}(y), G^{k}(z)), \\
  T^{3k+1}(x,y,z) &=& (f(F^{k}(y)), g(G^{k}(z)), h(H^{k}(x))), \\
  T^{3k+2}(x,y,z) &=& (f(g(G^{k}(z))), g(h(H^{k}(x))),h(f(F^{k}(y)))),
\end{eqnarray}
for $k \geq 0$, with $H^{0}, F^{0}$ and $G^{0}$ are the identity.\\
The cycles of the one-dimensional functions $ H$, $ F $ and $ G $ are linked; indeed, if $ x $ is a fixed point of $ H $, then $ y = g (h (x)) $ is a fixed point of $ F $ and $ z = h (x) $ is a fixed point of $ G $. Taking into account a cycle of order $ n $ of $ H$, the cycles of the functions $F $ and $ G $ will be called conjugate cycles, and are such that, if $ X = \{ x_{i} \}_{i = \overline{1, n}}$  is a cycle of order $n$ of $H$, then a conjugate cycle of $F$ exists, and is given by $Y = g(h (X ))$ and the conjugate cycle of $G$ is given by $Z=h(X)$.
\subsection{Homogeneous cycles}
\begin{definition}
A cycle of the map $T$ is said homogeneous if the components of its periodic points belong to conjugated cycles of $H$, $F$ and $G$. Otherwise, it is
called mixed cycle.
\end{definition}

\begin{Proposition}\label{prph301}\cite{A27}
If $n$ is not a multiple of $3$, the homogeneous cycles of period $n$ of $3D$ map $T$ associated with the cycle $ X = \{x_{1}, x_{2}, ..., x_{n} \} $ of $1D$ map $H$ of the same period have only one periodic point with first component $ x_{1}$. This periodic point is $(x_{1},y_{2s+1},z_{s+1})$  when  $n = 3s+1$; and $(x_{1},y_{s+1},z_{2s+2})$  when $n = 3s+2$.
\end{Proposition}
\begin{Proposition}\label{prph302}\cite{A27}
All the different homogeneous cycles of period $3n$ of the $3D$ map $T$ obtained starting from a cycle of period $n  \geq  2$ of the $1D$ map $H$ can be obtained by the periodic points:
         \begin{center}
         \begin{quote}
         $(x_{1},y_{j},z_{j+h})$  with  $h \leq j \leq n-2h $ \ and \ $ 1 \leq h \leq  \lfloor \frac{n}{3}   \rfloor $,\\
         $(x_{1},y_{j},z_{j+1-h})$  with  $2h-1 \leq j \leq n-h  $ \ and \ $ 1 \leq h \leq  \lfloor \frac{n+1}{3}   \rfloor $,
         \end{quote}
         \end{center}
where $\lfloor x   \rfloor$ is the ceiling function (the largest integer smaller than $x$).
\end{Proposition}

\subsection{Mixed cycles}
 Let $\{ x_{i}\}_{i=\overline{1,n}}$, $\{a_{j}\}_{j=\overline{1,m}}$ and $\{\alpha_{l}\}_{l=\overline{1,p}}$ three cycles of $H$ of period $n$, $m$ and $p$ (respectively), with their conjugated cycles ($y_{i}= g\circ h(x_{i})$, $b_{j}= g\circ h(a_{j})$ and $\beta_{l}= g\circ h(\alpha_{l})$ are the periodic points of the cycles of $F$, with $z_{i}= h(x_{i})$, $c_{j}= h(a_{j})$ and $\gamma_{l}= h(\alpha_{l})$). We  have:
\begin{Proposition}\label{prp303}\cite{A27}
If the map $H$ has two distinct coexisting cycles of period $n$ and $m$ and $s =LCM(n,m)$, then the map $T$ has $(n+m)\frac{n.m}{s}$ different mixed cycles of period $3s$, besides the homogeneous one. All the distinct mixed cycles can be obtained by the periodic points:
\begin{center}
$(x_{1},b_{j},z_{l})$ \ with \ $1\leq j \leq d$ \ and \ $ 1 \leq l \leq n $,\\
$(x_{1},b_{j},c_{l})$ \ with \ $1\leq j \leq d$ \ and \ $ 1 \leq l \leq m $.
\end{center}
\end{Proposition}
\begin{Proposition}\label{prp304}\cite{A27}
If the map $H$ has three coexisting cycles of period $n$, $m$ and $p$, and $S =LCM(n,m,p)$, then map $T$ has $2\frac{n.m.p}{S}$ different mixed cycles
of period $3S$, besides the homogeneous ones and those generated by any pair of the three cycles of $H$. All the distinct mixed cycles obtained by mixing of all the components of the three cycles can be obtained by the periodic points:
\begin{center}
$(x_{1},b_{j},\gamma_{l})$ \ with \ $1 \leq j \leq d$ \ and \ $ 1 \leq l \leq p\frac{LCM(m,n)}{S} $,\\
$(x_{1},\beta_{l},c_{j})$ \ with \ $1 \leq j \leq d$ \ and \ $ 1 \leq l \leq p\frac{LCM(m,n)}{S} $,
\end{center}
where $d=gcf(n,m)$.
\end{Proposition}
\section{Study of cycles of the $3D$ map $T$}
For $k = 1$, the construction of the cycles will be done as follows:
\begin{eqnarray*}
T^{3}(x,y,z) &=& (H(x), F(y), G(z)) = (x^{2}+b ,y^{2}+b,  z^{2}+b).\\
T^{4}(x,y,z)&=&(f(F(y)), g(G(z)), h(H(x))) =(y^{2}+b,z^{2}+b, (x^{2}+b)^{2}+b).\\
T^{5}(x,y,z) &=& (f(g(G(z))), g(h(H(x))),h(f(F(y))))  =  (z^{2}+b,(x^{2}+b)^{2}+b,(y^{2}+b)^{2}+b).
\end{eqnarray*}
\subsection{Fixed points of  $T$}
From the fixed points of the $1D$ map $H$ and its conjugate functions $F$ and $G$, we can define all homogeneous cycles of order 1  and all mixed cycles of order 3 of the $3D$ map $T$.
\begin{equation*}
 H(x) = f(g(h(x)))= x^{2}+b.
\end{equation*}
Fixed points of $T$ denoted $X_{1}$ and $X_{2}$ are:
\begin{eqnarray*}
  X_{1} &=& (x_{1},g(h(x_{1})),h(x_{1}))= (\frac{1}{2}+ \frac{1}{2}\sqrt{1-4b},\frac{1}{2}+ \frac{1}{2}\sqrt{1-4b},\frac{1}{2}+ \frac{1}{2}\sqrt{1-4b}), \\
  X_{2} &=& (x_{2},g(h(x_{2})),h(x_{2}))= (\frac{1}{2}- \frac{1}{2}\sqrt{1-4b},\frac{1}{2}- \frac{1}{2}\sqrt{1-4b},\frac{1}{2}- \frac{1}{2}\sqrt{1-4b});
\end{eqnarray*}
with $x_{1}=\frac{1}{2}+ \frac{1}{2}\sqrt{1-4b}$ for $ b \leq \frac{1}{4} $ is a fixed point of $H$, $y_{1}=g(h(x_{1}))$ and $z_{1}=h(x_{1})$ are conjugate fixed points of $F$ and $G$ respectively. Similarly, $x_{2}=\frac{1}{2}- \frac{1}{2}\sqrt{1-4b}$ pour $ b \leq \frac{1}{4} $ is a fixed point of $H$, $y_{2}=g(h(x_{2}))$ and $z_{2}=h(x_{2})$ are the conjugate fixed points of $F$ and $G$  respectively.\\
 Conjugate cycles have the same stability, because they have equal multipliers. The multipliers associated with the conjugated cycles are:
 \begin{equation*}
 \lambda _{x} = (H(x_i))^{'}= 2x_{i}, \  \lambda _{y} = (F(g(h(x_i))))^{'}= 2y_{i}\ \text{et}\  \lambda _{z} = (G(h(x_i)))^{'}= 2z_{i}, i=1, 2.
 \end{equation*}
 The multipliers associated with the conjugate cycles are equal to the eigenvalues of the Jacobian of $T^{3}$. Indeed, the Jacobian of $T^{3}$ is:
\begin{equation}
J_{T3}=\left(
     \begin{array}{ccc}
       2x & 0 & 0 \\
       0 & 2y & 0 \\
       0 & 0 & 2z
     \end{array}
\right)
\end{equation}
The characteristic polynomial $P_{1}(\lambda)= (2x-\lambda)(2y-\lambda)(2z-\lambda)$ vanishes for $\lambda = 2x$, $\lambda = 2y$ and $\lambda = 2z$ and
the stability of the fixed points of $H$ induces that of the conjugates of $F$ and $G$, because they have the same eigenvalues.
Likewise for the $3D$ maps $T$ and $T^{3}$, the stability of the fixed points of $ H $ induces that of the fixed points and third order cycles of $T$.
Knowing that the fixed points $x_{1}$  and $x_{2}$ of $H$ are unstable and stable respectively for $\frac{-3}{4} <  b < \frac{1}{4}$, then the fixed point $ X_{1}$ of $T$ generated by $x_{1}$ is unstable and $X_{2}$ generated by $x_{2}$ is stable under the same conditions on $b$.

\subsection{Bifurcation of fixed points of the $3D$ map $ T $}
Similarly, the $ T $ bifurcations are deducted from that of $ H $.

\indent {\bf $\bullet$ Fold Bifurcation:}  The Fold bifurcation of $H$ corresponds to  $ H(x)=x $ and   $ H'(x)=1$, then $ T $ undergoes a Fold bifurcation at  $b=\frac{1}{4}$ and $x = \frac{1}{2}$, $y = g(h(x))$ and $z = h(x)$.\\
\indent {\bf $\bullet$ Flip Bifurcation:}  The Flip bifurcation of $H$ corresponds to  $ H(x)=x $ and   $ H'(x)=-1$, then $T$ undergoes a Flip bifurcation at $b=\frac{-3}{4}$ and $x = y = z =\frac{-1}{2}$.\\
\indent {\bf $\bullet$ Transcritical bifurcation of fixed points of $T$:} A transcritical bifurcation of the fixed points  $x_{1}$ and $x_{2}$ of $H$ occurs when there exists a $\ b^{*}$ such that for $b > b^{*}$ the two fixed points exchange their stability, with $H^{'}(x_{1})= H^{'}(x_{2})=1$, for $b= b^{*}$. Here $T$ checks this property for $ b^{*} = \frac{1}{4}$.

\subsection{Cycles of period 2}
From the proposition \ref{prph301}, the map $ T $ admits a homogeneous cycle of period two $C_{2}=\{V_{1}, V_{2}\}$, deduced from the cycle of order two $c_{2}=\{V_{1x}, V_{2x}\}$ of $H$, where $V_{1x} = \frac{-1}{2}+\frac{1}{2}\sqrt{-3-4b}$ and $V_{2x} = \frac{-1}{2}-\frac{1}{2}\sqrt{-3-4b}$. The homogeneous cycle of period two of $T$, whose components are $ V_{1} =(V_{1x},V_{2x},V_{1x})$ and $V_{2} = (V_{2x},V_{1x},V_{2x})$, is stable for $\frac{-5}{4}< b < \frac{-3}{4}$.\\
The Flip bifurcation of the two-period cycle of $T$ is deduced from the Flip bifurcation of the two-period cycle of $H$ ($H^{2}(x)=x$ and $(H^{2})^{'}=-1$): $b = \frac{-5}{4}$ for $x \in \{\frac{1}{2}\sqrt{2}-\frac{1}{2}, -\frac{1}{2}\sqrt{2}-\frac{1}{2}\}$.\\
The  cycle of the period two mixed with an fixed point of $1D$ map $ H $  builds a  cycle of period six  of the $3D$ map $T$.

\subsection{Construction of the period cycles greater than 2 of $T$}
 Using the cycles of the $1D$ map $H$ and the conjugate cycles of $F$ and $G$, the homogeneous and mixed cycles of the $3D$ map $T$ are constructed.

\subsubsection{Cycles of periods 4 and 5}
\begin{figure}[!ht]
\begin{center}
    \subfigure[Cycles of periods 1 to 4]{\label{H4}\includegraphics[width=3.5cm,height=3.3cm]{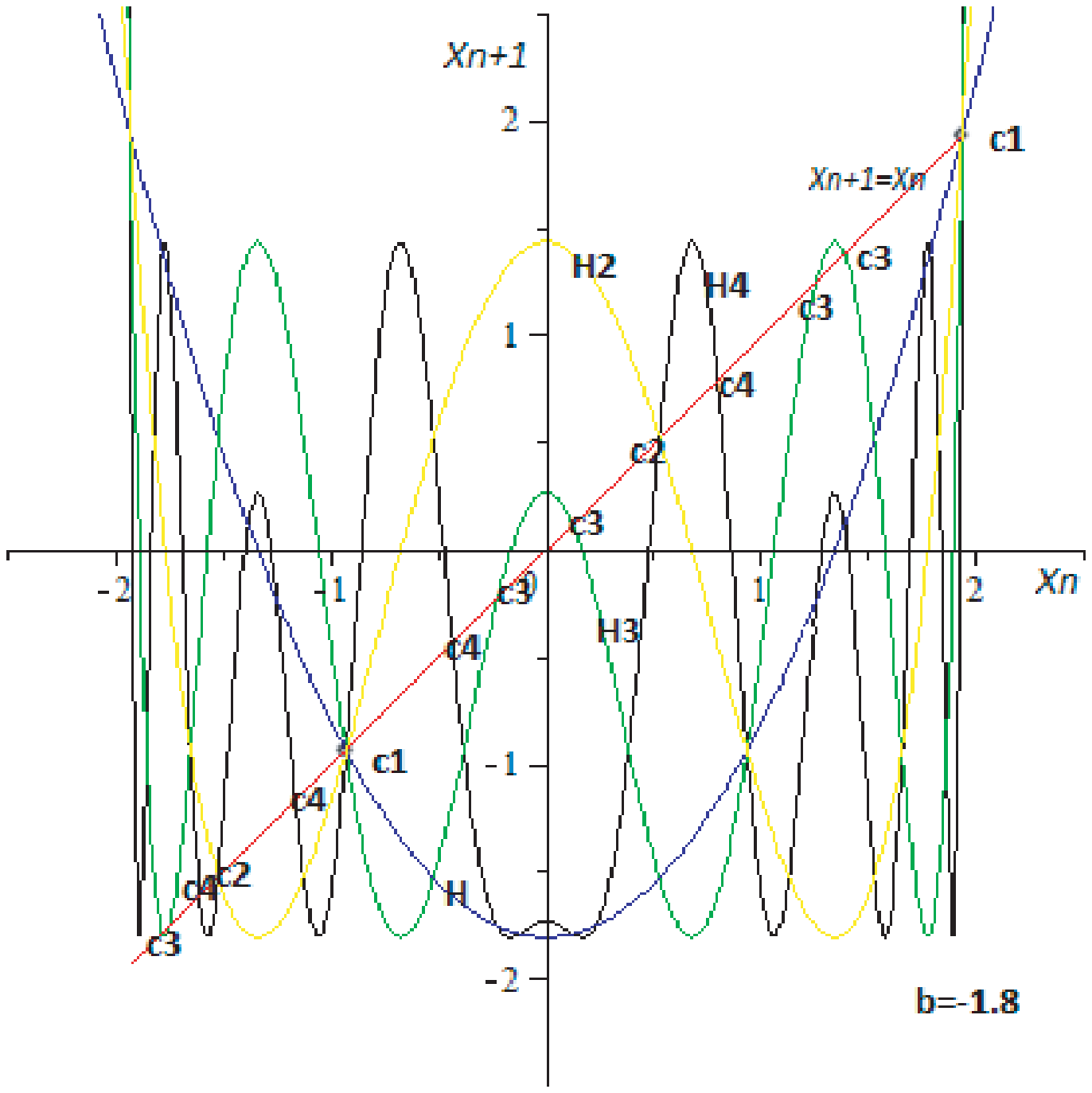}}
    \hspace{0.5cm}
    \subfigure[Stability]{\label{C1a4}\includegraphics[width=3.7cm,height=3.5cm]{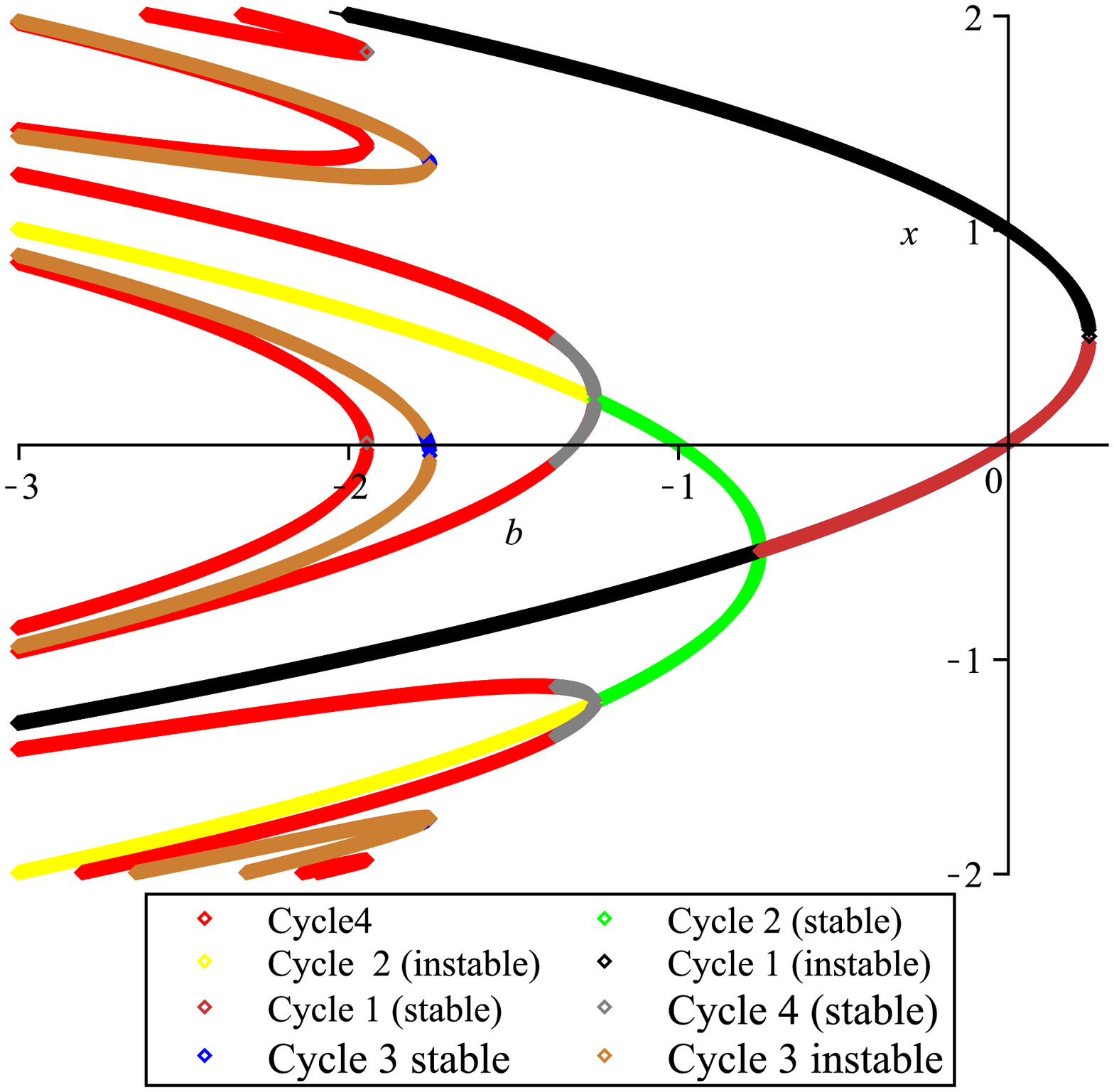}}
    \hspace{0.5cm}
     \subfigure[Cycle of period 6]{\label{C6xyzn-a}\includegraphics[width=3.5cm,height=3.3cm]{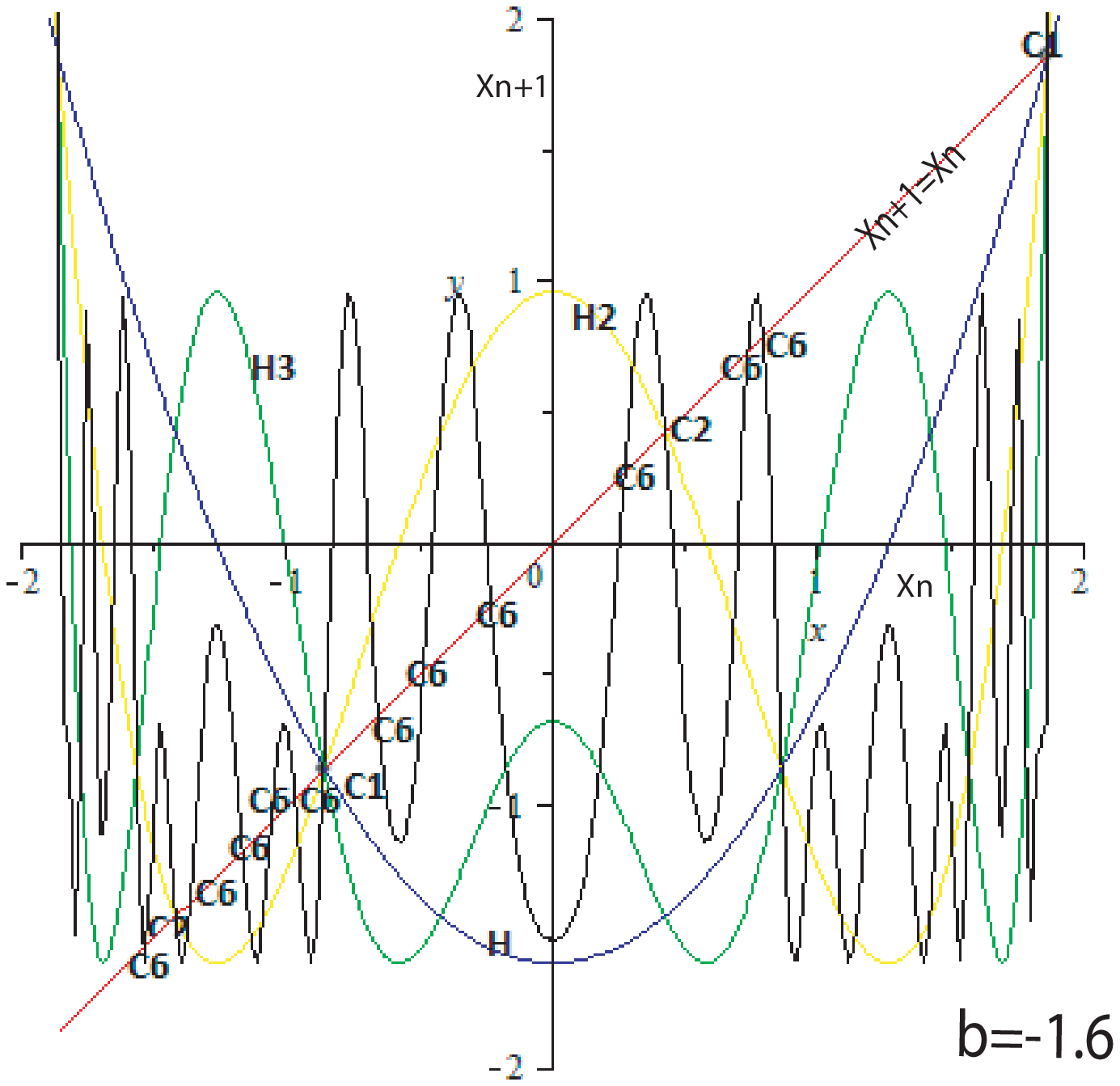}}
    \vspace{-0.3cm}
     \caption{A few cycles of $ H $ and their stability.}\label{H4n}
 \end{center}
\end{figure}

Under the proposition \ref{prph301}, the transformation $T$ has a unique homogeneous cycle of order four $C4$, and its components are derived from the  cycle of the period four of $H$, $c4= \{ a_{1},a_{2},a_{3},a_{4} \}$ (points noted $c4$ in the figure \ref{H4}), using the periodic point $(a_{1},a_{4},a_{3})$. The cycle of period four of $T$ is stable for $-1.38 < b < \frac{-5}{4}$ (figure \ref{C1a4}).\\
Similarly, there is a single homogeneous cycle $ C5 $ of period five of  and its components are deduced from the cycle of period five of $H$,  $c5=\{ b_{1},b_{2},b_{3},b_{4}, b_{5} \}$. The cycle components of five $C5$ periods of $T$ are constructed using the periodic point $(b_{1},b_{3},b_{5})$.

\subsubsection{Cycles of period 3}
By virtue of the proposition \ref{prph301}, the point transformation $T$ admits only two mixed cycles of period $3$, deduced from the fixed points of $H$, $x_{1}$ and $x_{2}$. The cycles of  period three of $T$, are given by:
\begin{description}
  \item[$\bullet$] $C3_{1} =\{V3_1,V3_2,V3_3\}$:
  \begin{eqnarray*}
  V3_1 &=& (x_{2},x_{1},x_{2})= [\frac{1}{2}-\frac{1}{2}\sqrt{1-4b}, \frac{1}{2}-\frac{1}{2}\sqrt{1-4b}, \frac{1}{2}+\frac{1}{2}\sqrt{1-4b}],  \\
  V3_2 &=& T(V3_1) \ \text{et} \ V3_3=T(V3_2).
\end{eqnarray*}
  \item[$\bullet$] $C3_{2}=\{V3_4,V3_5,V3_6\}$:
\begin{eqnarray*}
  V3_4 &=& (x_{2},x_{1},x_{1})=[\frac{1}{2}-\frac{1}{2}\sqrt{1-4b}, \frac{1}{2}+\frac{1}{2}\sqrt{1-4b}, \frac{1}{2}+\frac{1}{2}\sqrt{1-4b}], \\
  V3_5 &=& T(V3_4) \ \text{et} \ V3_6=T(V3_5).
\end{eqnarray*}
\end{description}
The Flip bifurcation of the three-period cycle of $T$, occurs for $b= -1.768529152$; and the Fold bifurcation of a three-period cycle of $T$ occurs for $b = \frac{-7}{4}$.

\subsubsection{Cycles of period 6}
Under the stated propositions, there is a total of nine cycles of period six for the map $T$; indeed:
 \begin{enumerate}
\item From the proposition (\ref{prph302}), a six-period homogeneous cycle of $ T $ deducted from the two-period cycle $c2=\{ V_{1x}, V_{2x}\}$ of $H$, is obtained using the periodic point $(\alpha_{1},\alpha_{1},\alpha_{2})$, with $\alpha_{1}=V_{1x}$ and $\alpha_{2}=V_{2x}$; stable for $\frac{-5}{4} < b < \frac{-3}{4}$.
\item From the proposition (\ref{prp304}), we have the existence of eight mixed cycles of period six of $T$:
\begin{enumerate}
 \item There are two mixed cycles of six period deduced from the coexistence of the two fixed points, $\{x_{1}\}$ and $\{x_{2}\}$
  and the period cycle two $c2$ of $H$ (here we  have $n=1$, $m=1$ and $p=2$, then $S = p.p.c.m(n,m,p)= 2$). The first component of each of the two cycles of period  six is $(\alpha_{1},x_{2},x_{1})$ and $(\alpha_{1},x_{1},x_{2})$ respectively.
\item From the coexistence of the cycles $\{x_{1}\}$ and $c2$ of $H$, three periodic cycles of period six are obtained from the periodic points: $(\alpha_{1},x_{1},\alpha_{1})$, $(\alpha_{1},x_{1},\alpha_{2})$ and $(\alpha_{1},x_{1},x_{1})$.
\item And from the coexistence of the cycles $\{x_{2}\}$ et $c2$ de $H$, three other periodic cycles of period six are obtained from the periodic points: $(\alpha_{1},x_{2},\alpha_{1})$, $(\alpha_{1},x_{2},\alpha_{2})$ and $(\alpha_{1},x_{2},x_{2})$.
\end{enumerate}
\end{enumerate}
Note that it is still possible for us to construct mixed cycles with a period greater than six by using the proposition (\ref{prp304}).

\subsection{Bifurcation diagramme of the $3D$ map $T$}
\begin{figure}[!ht]
\begin{center}
    \subfigure[Bifurcation diagramme]{\label{cascade3-a}\includegraphics[width=5cm,height=3cm]{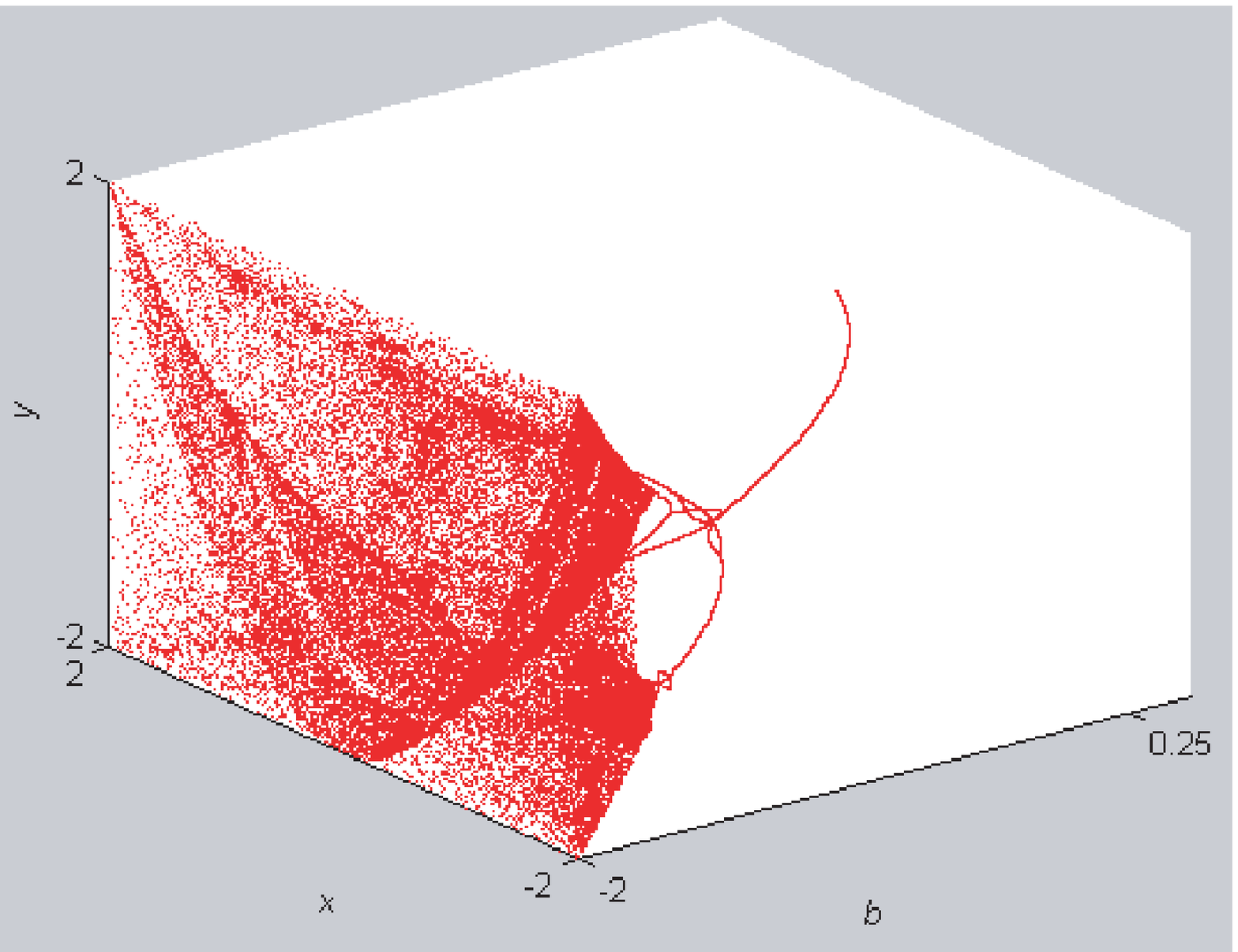}}
    \subfigure[Bifurcation diagramme in the plan $(b, x)$]{\label{cascade3-b}\includegraphics[width=5cm,height=3cm]{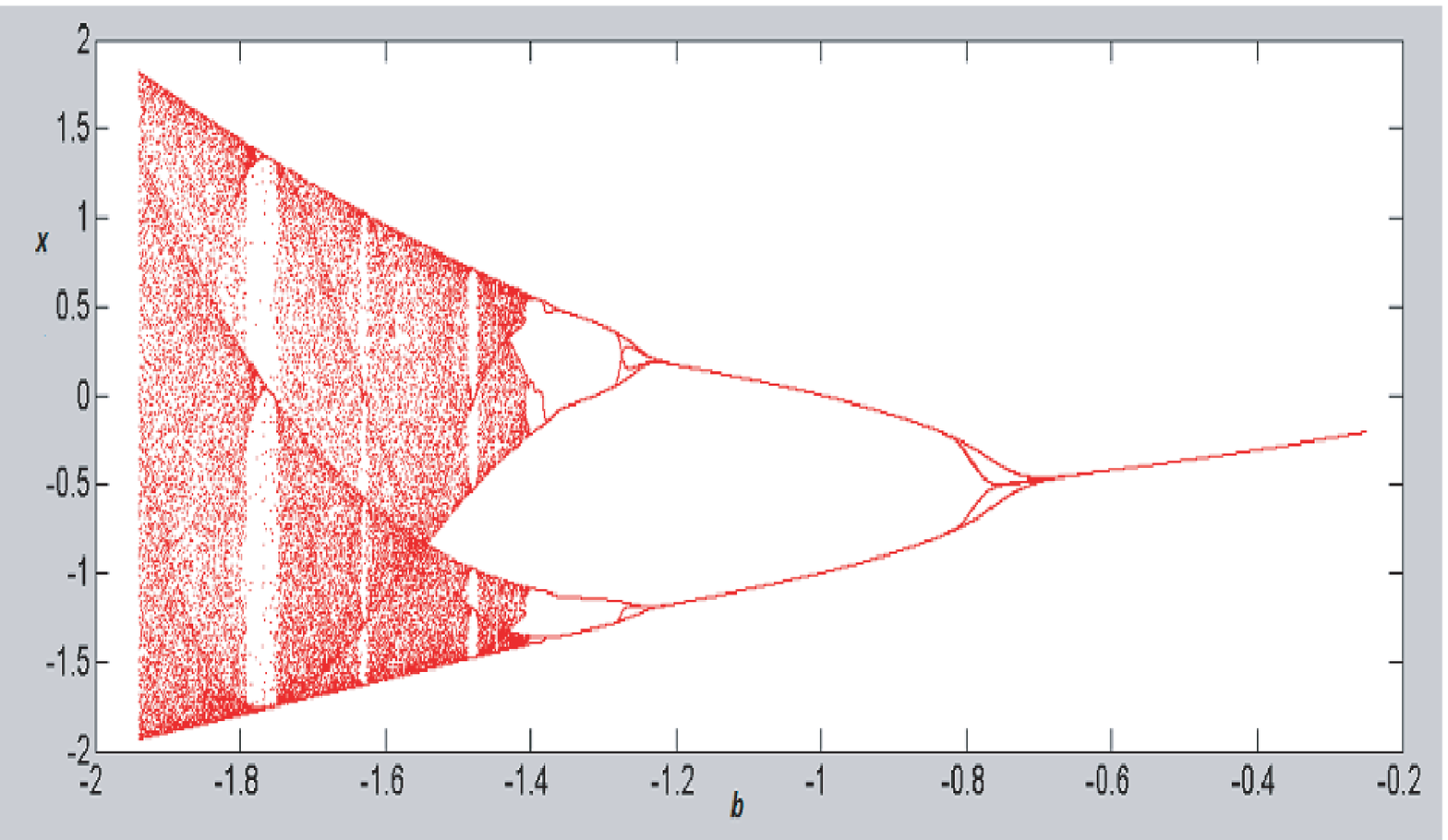}}
    \vspace{-0.3cm}
        \caption{Bifurcation diagramme of $T$ ($b \in [ -1.99, -0.3]$).}
  \label{cascade3}
\end{center}
\end{figure}

The figures \ref{cascade3-a} and \ref{cascade3-b} give an example of a Feigenbaum bifurcation diagram of the transformation $T$ defined in the space of the variables $x$, $y$ and $b$, and its projection in the plane $[b, x]$. As shown in the figure \ref{C1a4}, when the parameter $ b $ varies, we can observe the bifurcation in an attractive fixed point in an attractive  cycle of order 2, and then in an attractive  cycle order 4, etc.\\
In the figure \ref{cascade3-a},, we can see the projection of the bifurcation diagram in the plane $ [x, y] $, which represents the set of attractors of $ T $ (we chose here the initial value $X_{0} = (0, -0.5, 0 $).  Notice that the point $ (- 0.5, -0.5, -0.5) $ for $b=\frac{-3}{4}$, corresponds to a point in the doubly period bifurcation. The coexistence of the fixed points and the cycle of order two of $H$ gives rise to a cycle of order six of the $3D$ map $T$. From the bifurcation diagram (figure \ref{cascade3-b}), we can find the observations: evolution towards a fixed point for $ b = -0.4 $, a cycle of order two for $-0.8 < b < -0.75$, a cycle order four for $ b = -1.25 $. For $ b = -1.6 $, we no longer distinguish cycles; the system presents a chaotic character.

\section{Phase space of the $3D$ map $T$}
In this section we will give some attractors of the three-dimensional map $T$, delimited by critical manifolds and an example of the chaotic dynamics of $T$ by varying $b$.
\subsection{Critical manifolds}
\begin{figure}[!h]
  \begin{center}
  \includegraphics[width=7cm,height=3.5cm]{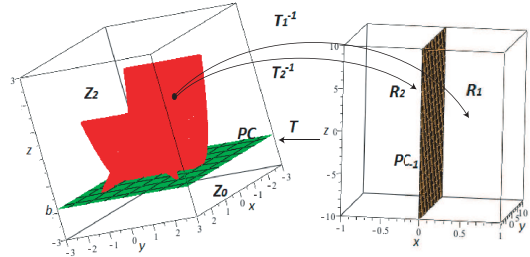}
  \end{center}
  \vspace{-0.5cm}
\caption{The regions $R_{1}$ and $R_{2} $ after applying inverse determinations of $T$.}\label{figF}
\end{figure}

\begin{figure}[!ht]
\begin{center}
    \subfigure{\label{fig:PC12-a}\includegraphics[width=3.5cm,height=3.5cm]{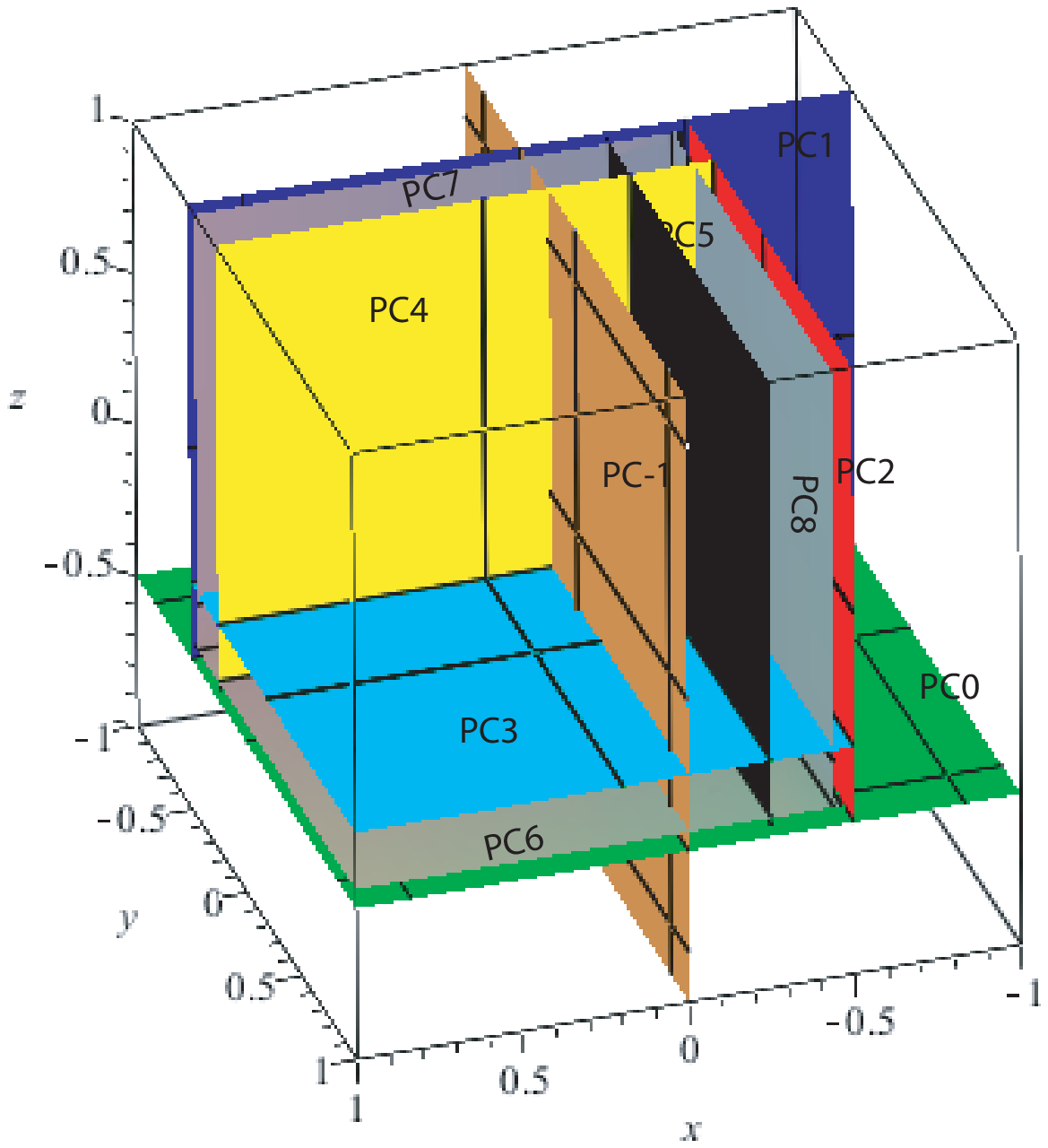}}
    \subfigure{\label{fig:PC12-b}\includegraphics[width=3.5cm,height=3.5cm]{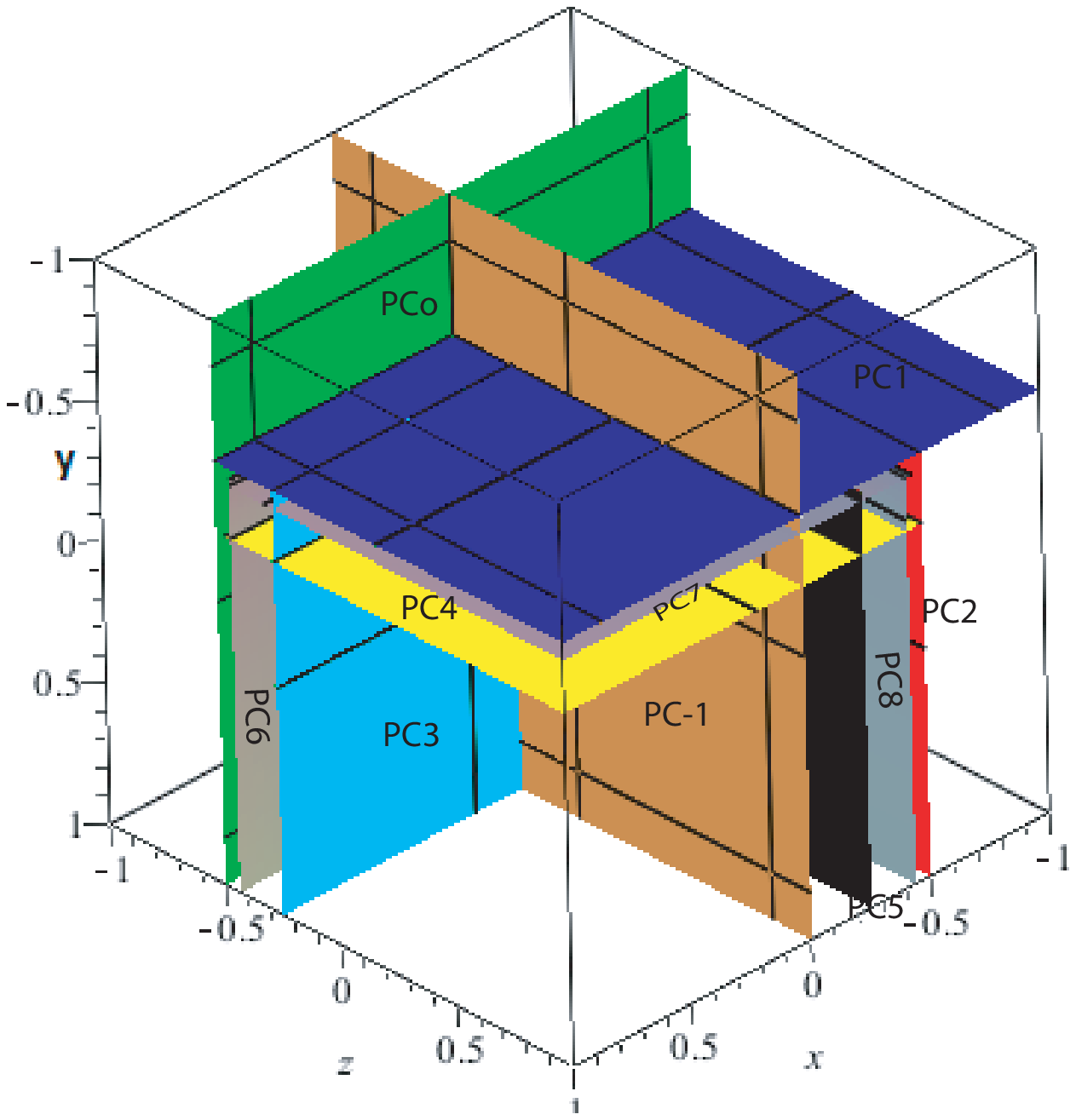}}
             \caption{A few critical plans of $T$ at $b=-0.5$.}
  \label{fig:PC12}
\end{center}
\end{figure}

We notice in the case of this map, that the critical manifolds are plans, denoted $PC$, dividing the space of the phases into zones denoted $ Z_{i}$, $i$ integer and each $Z_{i}$ is the set of points in the phase space that has $i$  antecedents of rank one. It is the generalization to dimension three of the notions of critical point and critical line defined in dimension one and two. A critical plan of rank $k+1$ noted $PC_{k}$ is the consequent plan of rank $k$ of $PC_{0}$,
$k=1,2,...$. $PC_{-1}$ is the antecedent of rank one of $PC_{0}$. These critical planes delimit the chaotic attractors.
The equation of the critical manifold $PC_{-1}$  of $T$ satisfies $\mid J_{1}(x,y,z)\mid = 0$, where $J(x,y,z)$ is the Jacobian of $ T $ at the point $(x,y,z)$.
\begin{equation}
J(x,y,z)=\left(
     \begin{array}{ccc}
       0 & 1 & 0 \\
       0 & 0 & 1 \\
       2x & 0 & 0
     \end{array}
\right)
\end{equation}
This critical variety is the plane $PC_{-1}=\{(0,y,z), y, z \in \mathbb{R} \}$.\\
The critical plane $PC=PC_{0}$ is therefore $ PC_{0} = T(PC_{-1})= \{(y, z, b),\ y, z \in \mathbb{R} \},$  it is a plane parallel to the plane $[x,y]$ which cuts the axis of $z$ into $b$; $PC$ separates the phase space into two regions. A region $Z_{2}$ verifying $z - b > 0$ and consisting of all points that have two rank $1$ antecedents; and a region denoted $Z_{0}$ such that $z -b < 0$ and whose points do not have antecedents.
We denote by $R_{1}$ and $R_{2}$ the regions verifying respectively $x>0$ and $x<0$, located on both sides of $PC_{-1}$ (figure \ref{figF}). The two inverse determinations corresponding to the two antecedents of rank $1$ of the points of the region $Z_{2}$ are:

\begin{eqnarray*}
  T^{-1}_{1}: Z_{2} \rightarrow R_{1}  \hspace{1.63cm}&    & T^{-1}_{2}: Z_{2} \rightarrow R_{2} \\
  T^{-1}_{1}:
\left\{
  \begin{array}{ll}
    x = \sqrt{z'-b}, & \hbox{} \\
    y = x', & \hbox{} \\
    z = y' .& \hbox{}
  \end{array}
\right. &    & T^{-1}_{2}:
\left\{
  \begin{array}{ll}
    x = -\sqrt{z'-b}, & \hbox{} \\
    y = x', & \hbox{} \\
    z = y' .& \hbox{}
  \end{array}
\right.
\end{eqnarray*}
Critical plans of order $n+1$, are defined by $PC_{n+1}=T(PC_{n})$  for all $n \geq 0$.
Equivalently for an order $n$, $ n \geq 0$, are defined by $PC_{n}=T^{n+1}(PC_{-1})$.
All critical plans for $T$, depend on the $b$ parameter, except $PC_{-1}$. The critical planes of $T$ are all parallel to one of the plane of the three-dimensional space $(x,y,z)$, as shown in figure \ref{fig:PC12}. The distance between two parallels critical planes  depend for  the value of $b$. For $b=0$, we can have critical plans confused, in this case, $PC_{-1}$, $PC_{2}$, $PC_{5}$ and $PC_{8}$.

\subsection{Chaotic Attractors}
By choosing the initial value $X_{0}=(0,0,-0.5)$, and $b=\frac{-3}{4}$, we have the appearance of an attractor formed of six branches that intersect at the point $(-0.5,-0.5,-0.5)$. At $b=-0.8$, the six branches split into two groups of three branches each, which forms an attractor of order six.
By decreasing $b$ up to $b = -1.45$, we have the appearance of six surfaces forming an attractor of order six. When $b$ reaches the value of $-1,864$, the attractor of order 6 becomes an attractor of order 3, consisting of a surface and two secant planes. For $b = -2$, the attractor is now composed of three line segments.

 \begin{figure}[!ht]
\begin{center}
    \subfigure[$b=-3/4$]{\label{AT1-a}\includegraphics[width=3.5cm,height=3.5cm]{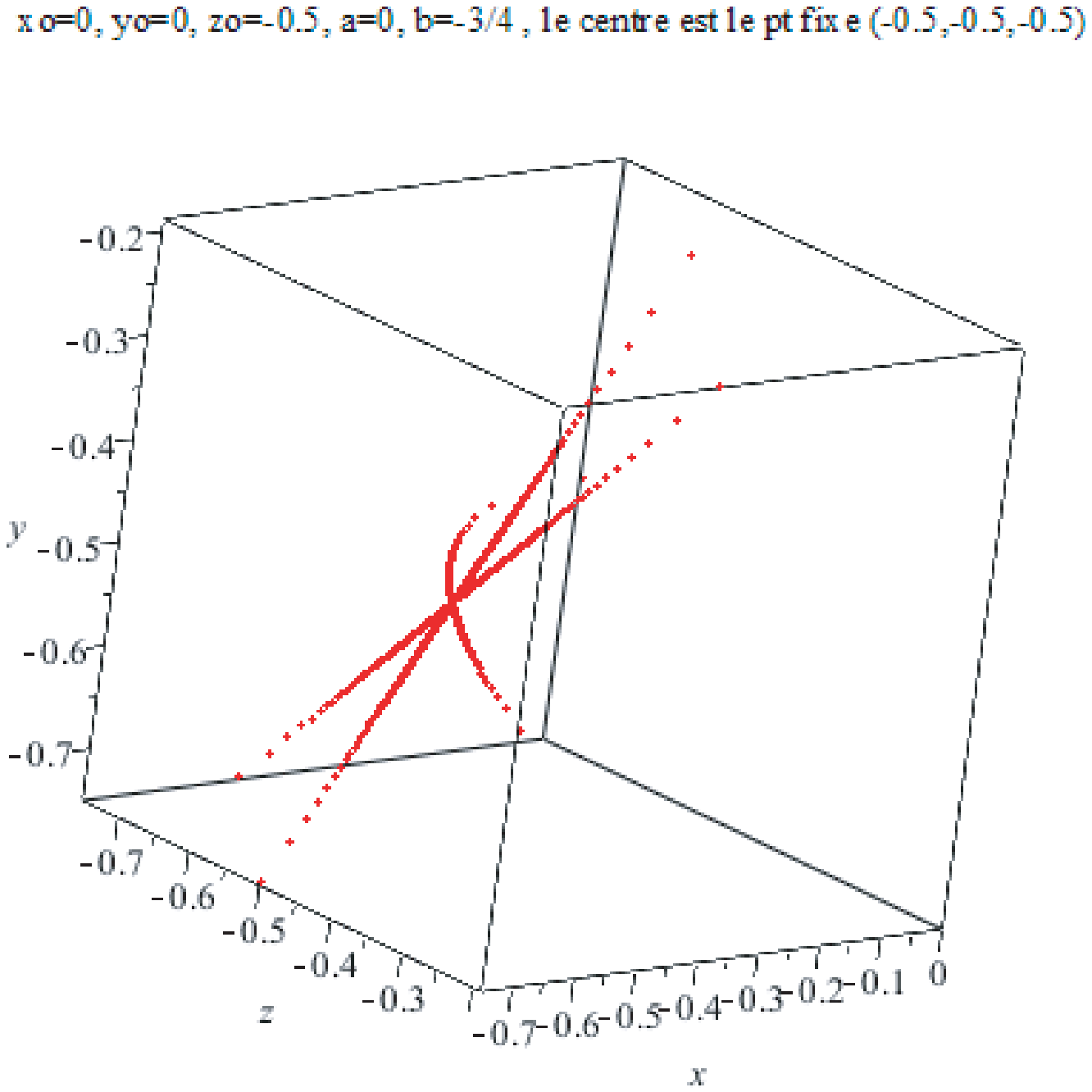}}
    \subfigure[$b=-0.8$]{\label{AT1-b}\includegraphics[width=3.5cm,height=3.5cm]{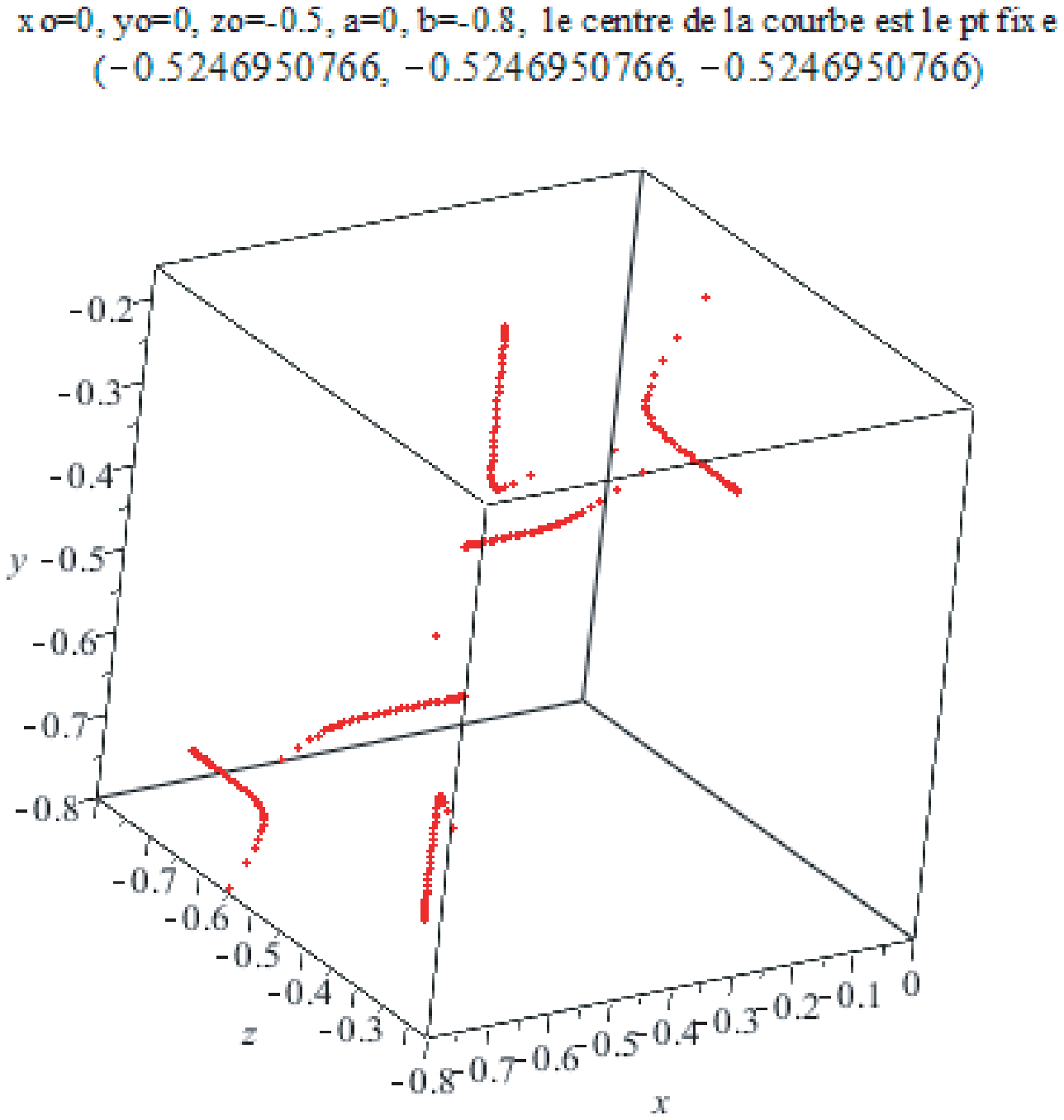}}
    \subfigure[$b=-1.45$]{\label{AT1-c}\includegraphics[width=3.5cm,height=3.5cm]{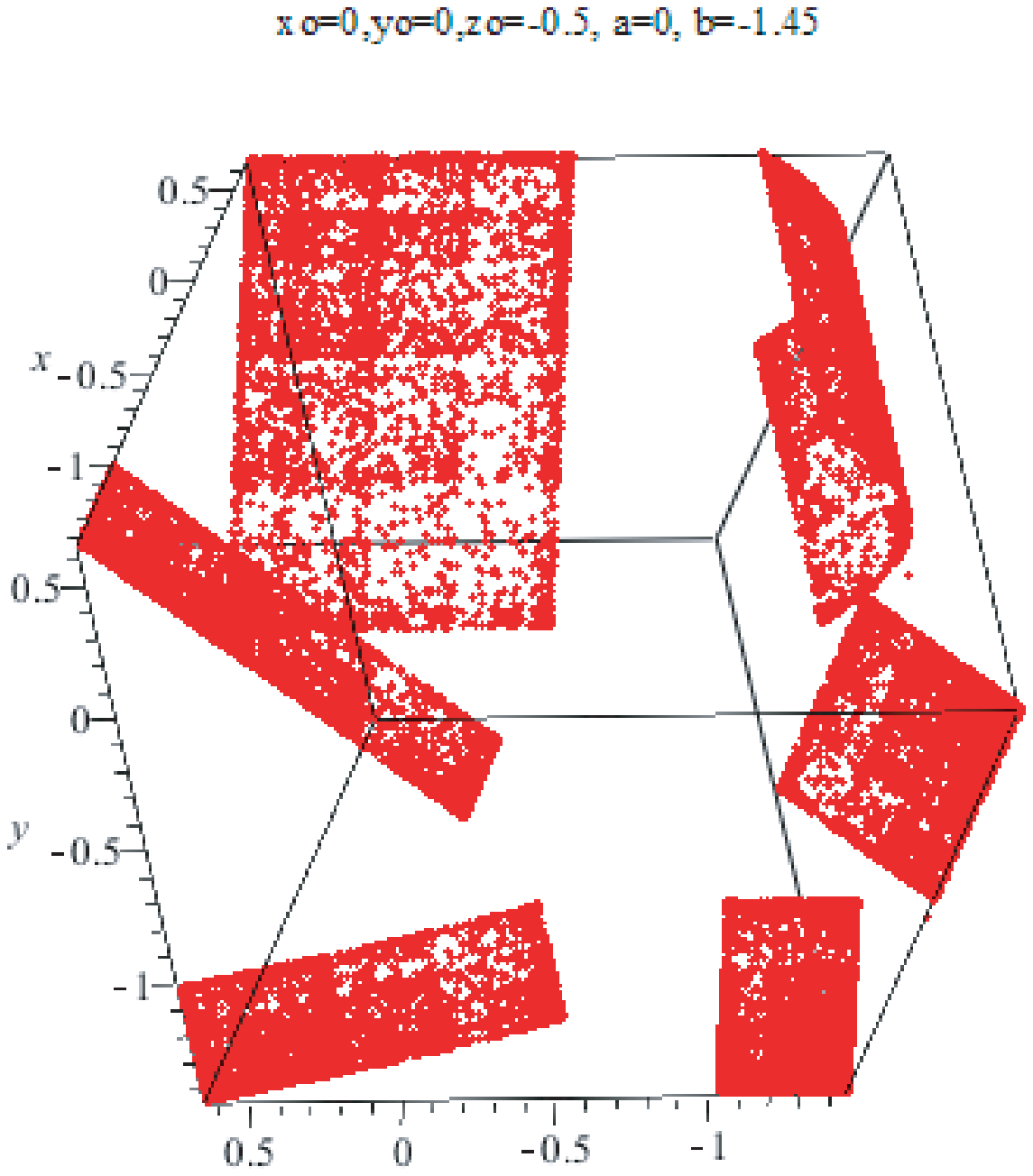}}\\
    \subfigure[$b=-1.54$]{\label{AT1-d}\includegraphics[width=3.5cm,height=3.5cm]{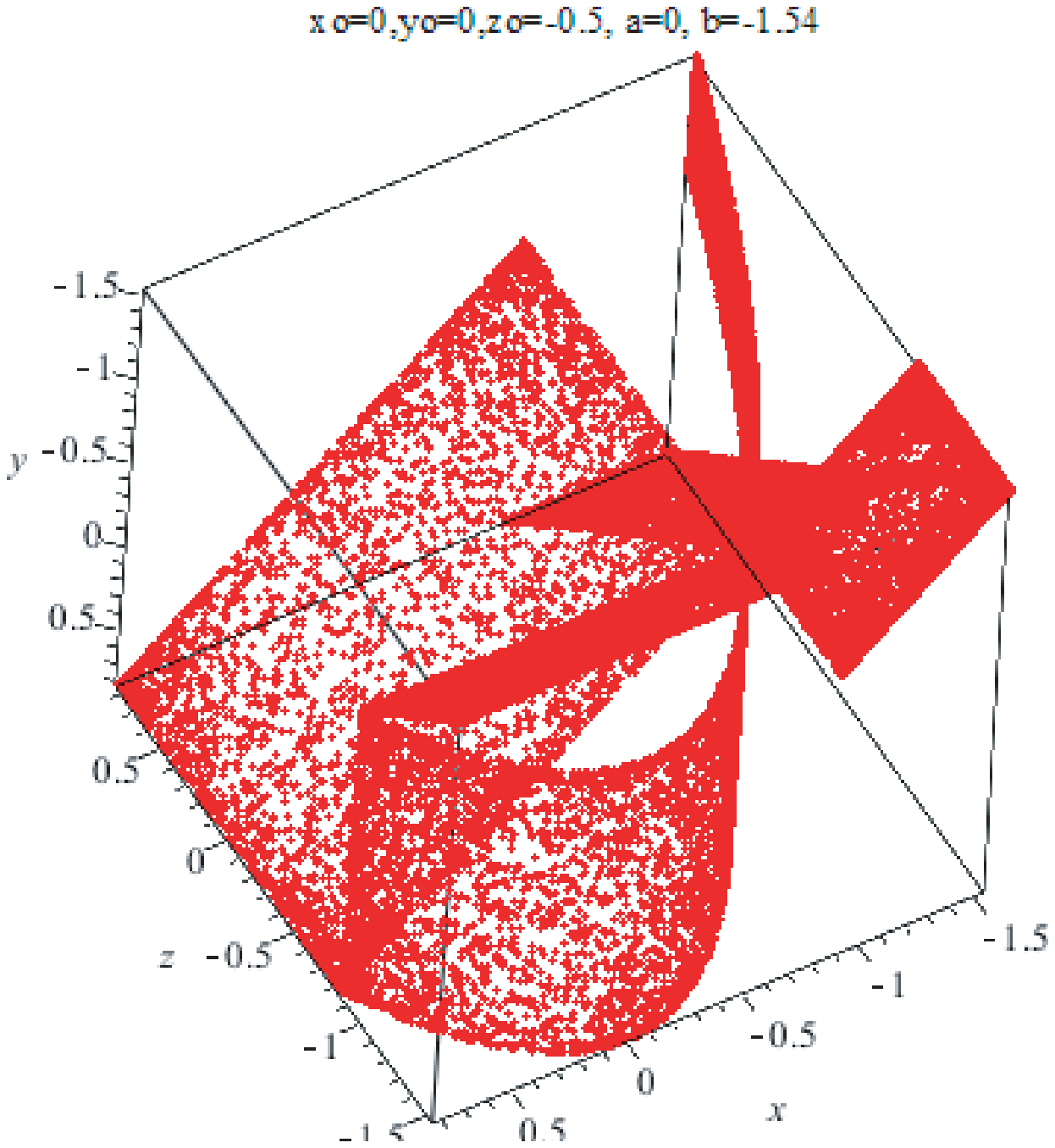}}
    \subfigure[$b=-1.864$]{\label{AT1-e}\includegraphics[width=3.5cm,height=3.5cm]{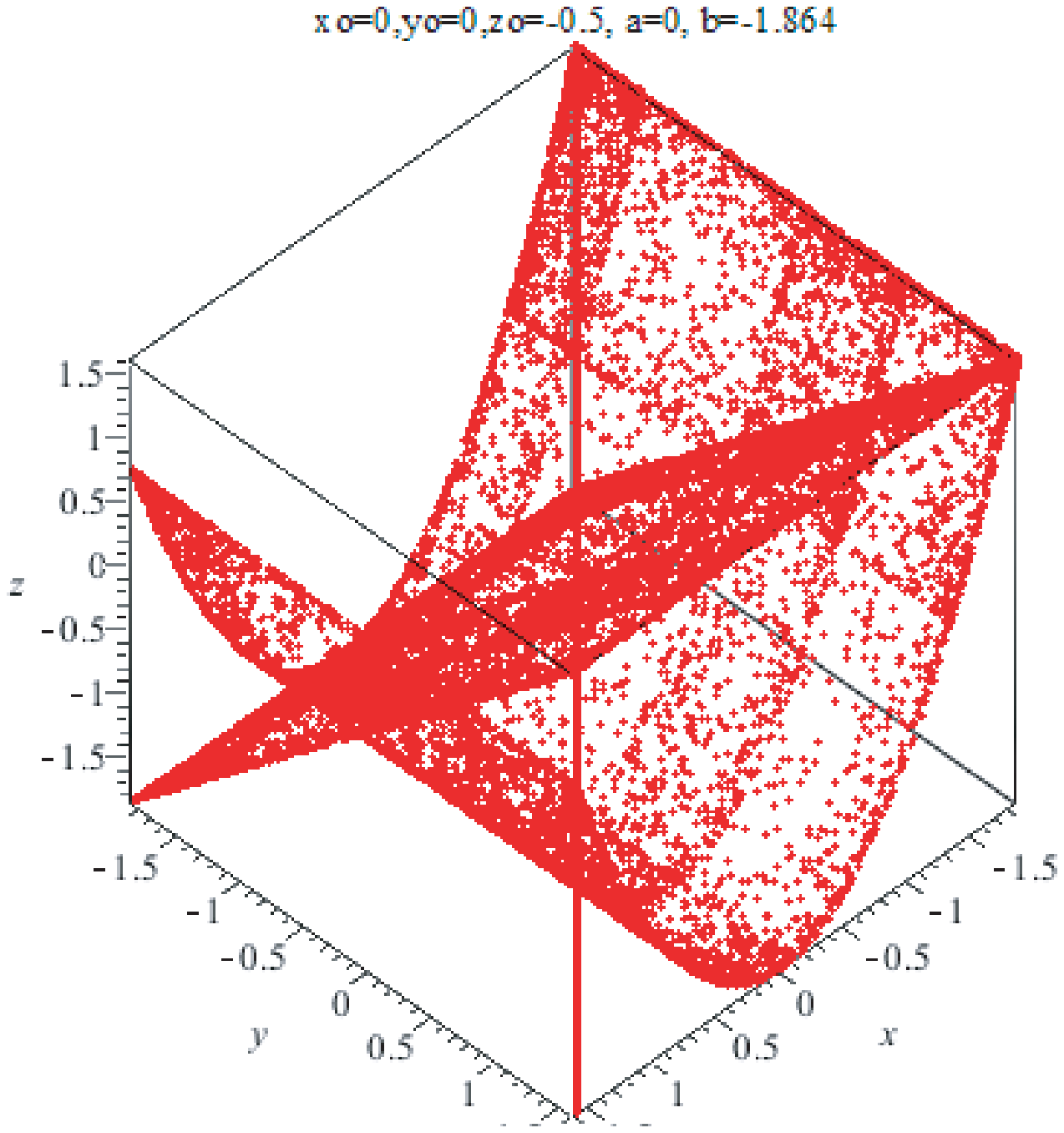}}
    \subfigure[$b=-2$]{\label{AT1-f}\includegraphics[width=3.5cm,height=3.5cm]{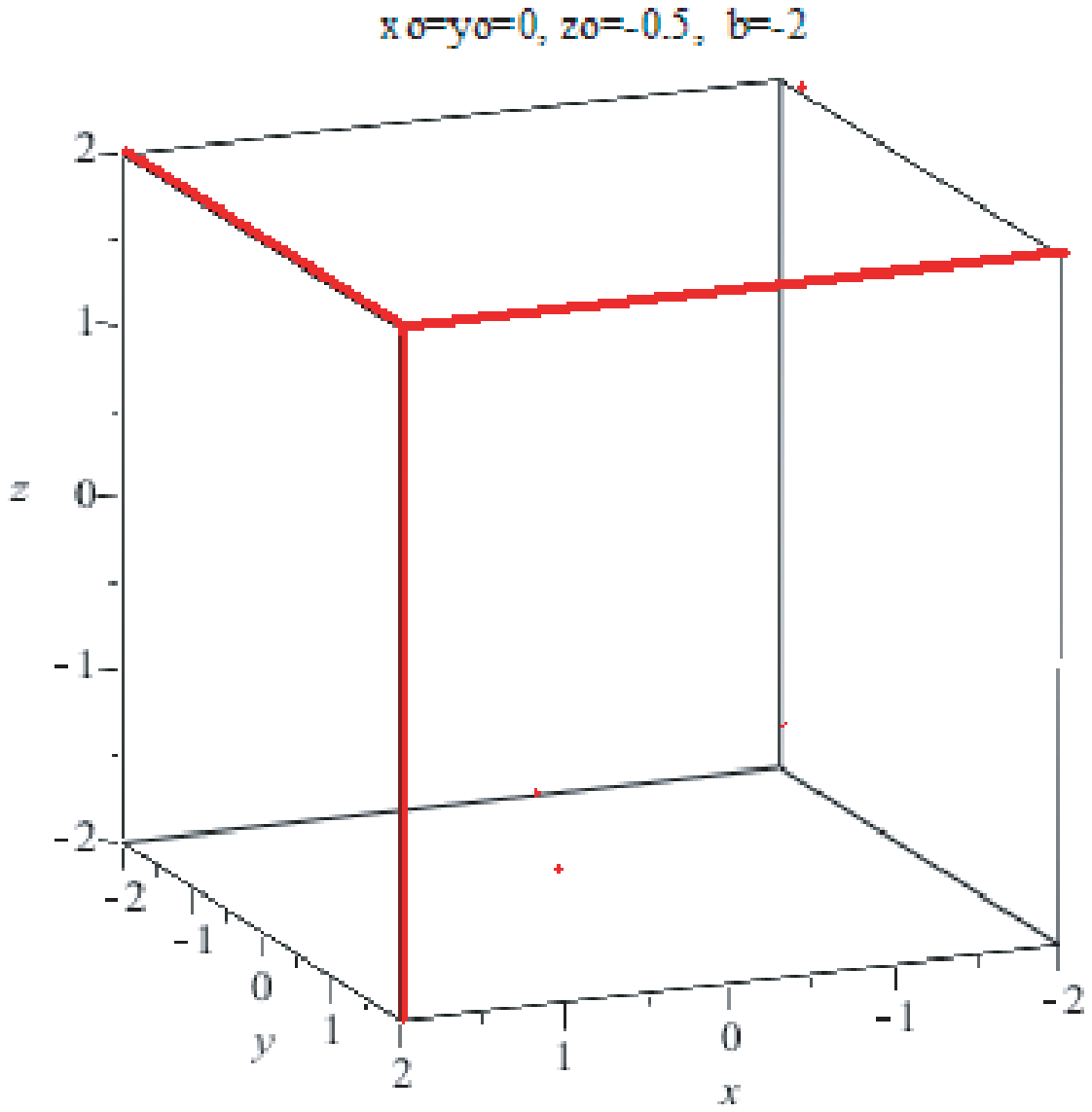}}
          \caption{Evolution of an attractor of order six to a attractor of order three.} \label{AT1}
\end{center}
\end{figure}

\begin{figure}[!ht]
\begin{center}
    \subfigure[$X_{0}(-0.5, -0.5, -0.48)$]{\label{At4-b}\includegraphics[width=3.5cm,height=3.5cm]{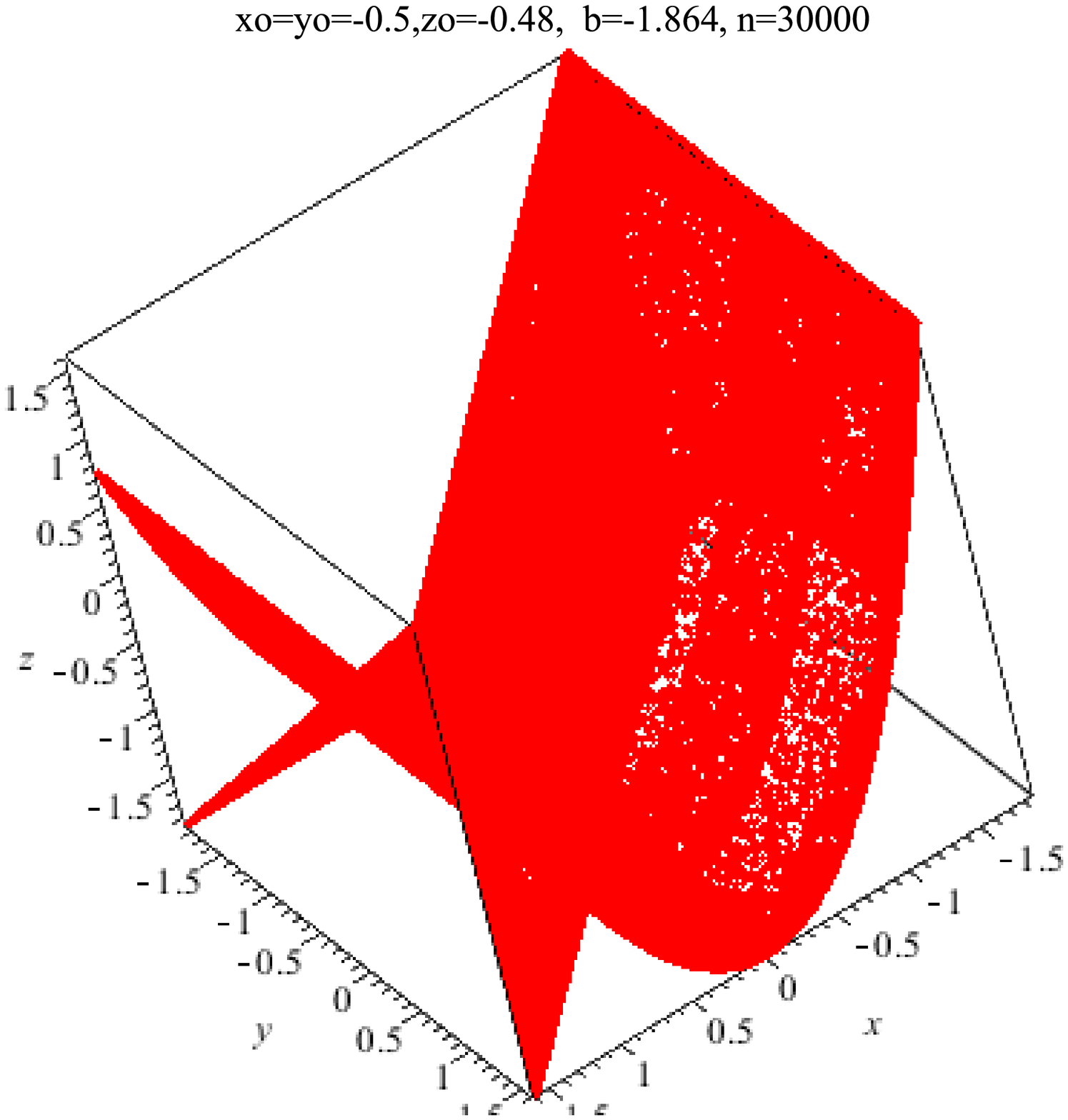}}
    \subfigure[$X_{0}(-0.5, -0.5, -0.5)$]{\label{At4-a}\includegraphics[width=3.5cm,height=3.5cm]{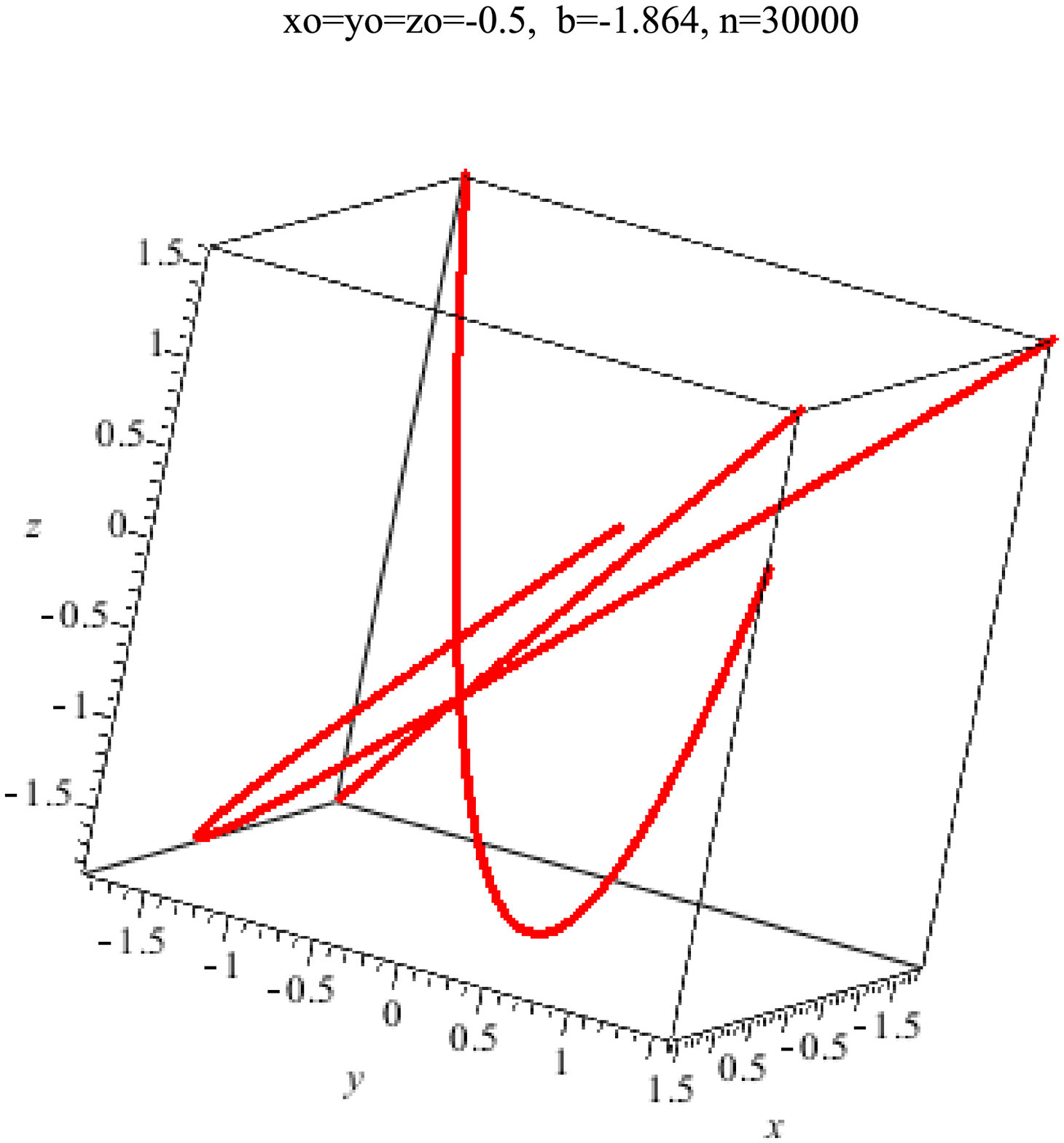}}
        \caption{The attractors of order six at  $b=-1.864$.} \label{At4}
\end{center}
\end{figure}

\begin{figure}[!ht]
\begin{center}
    \subfigure[$X_{0}(-0.5, 0, 0)$]{\label{At5-a}\includegraphics[width=3.5cm,height=3.5cm]{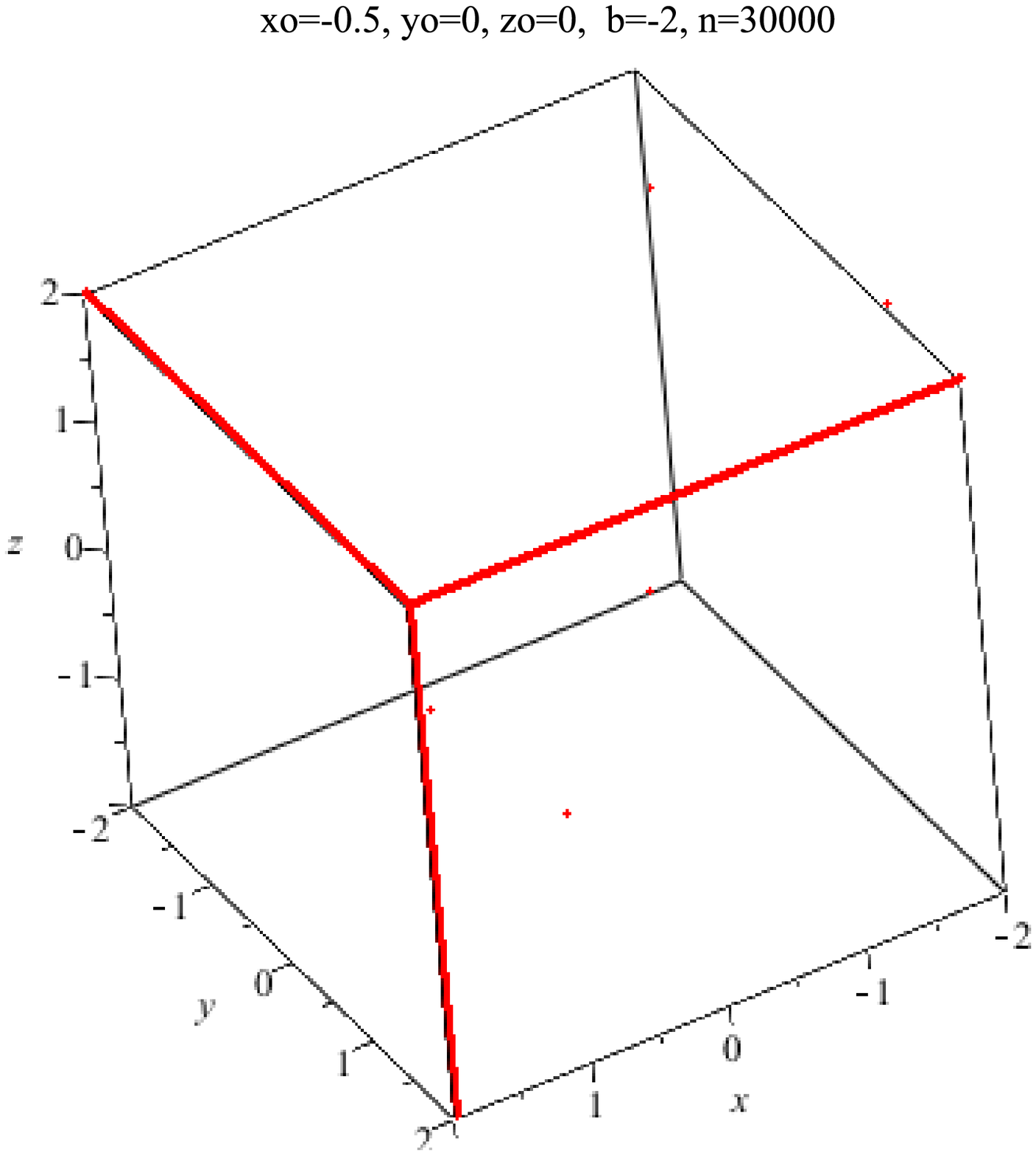}}
    \subfigure[$X_{0}(-0.5, -0.01, 0)$]{\label{At5-b}\includegraphics[width=3.5cm,height=3.5cm]{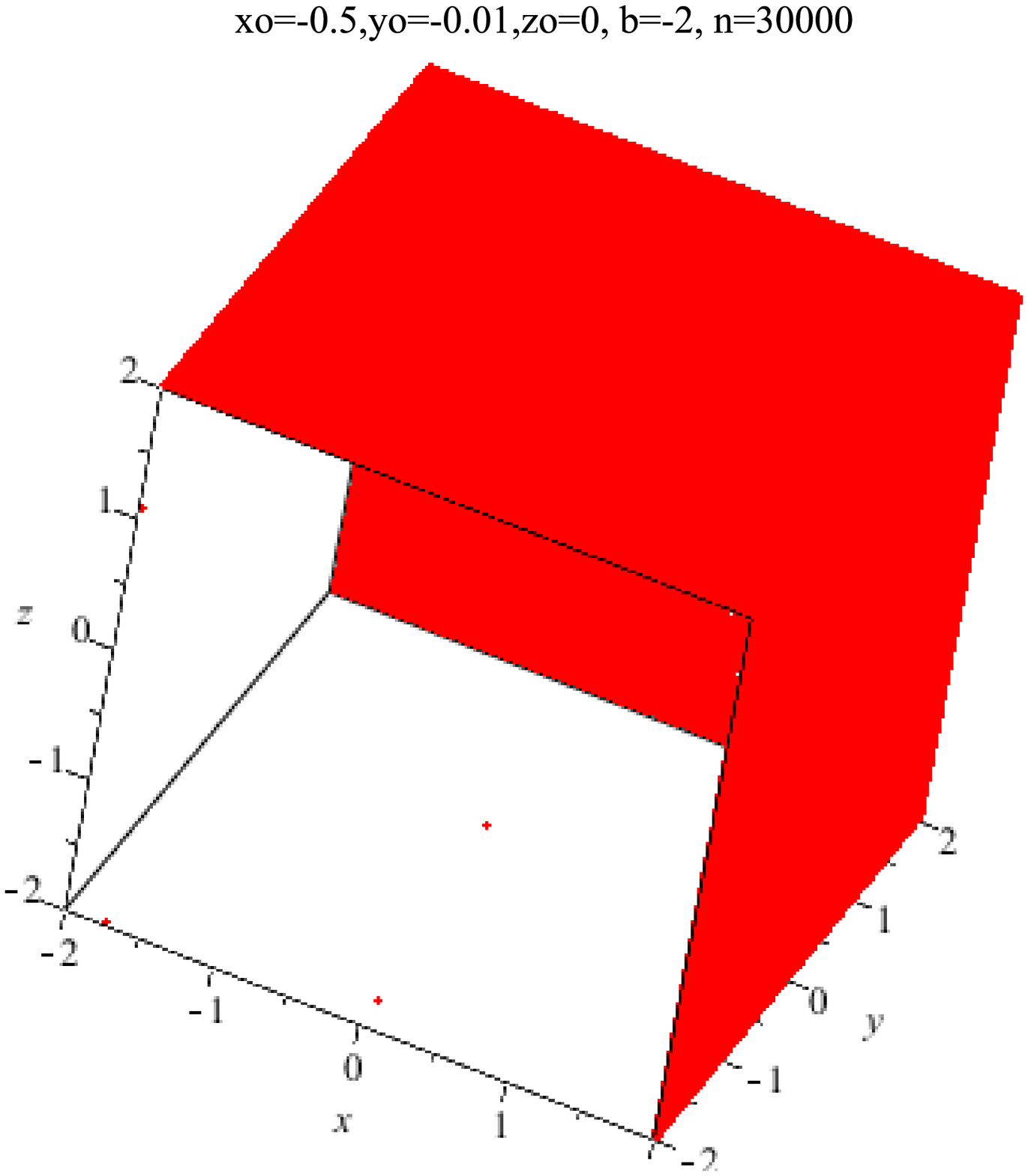}}
    \subfigure[$X_{0}(-0.5, -0.5, 0)$]{\label{At5-c}\includegraphics[width=3.5cm,height=3.5cm]{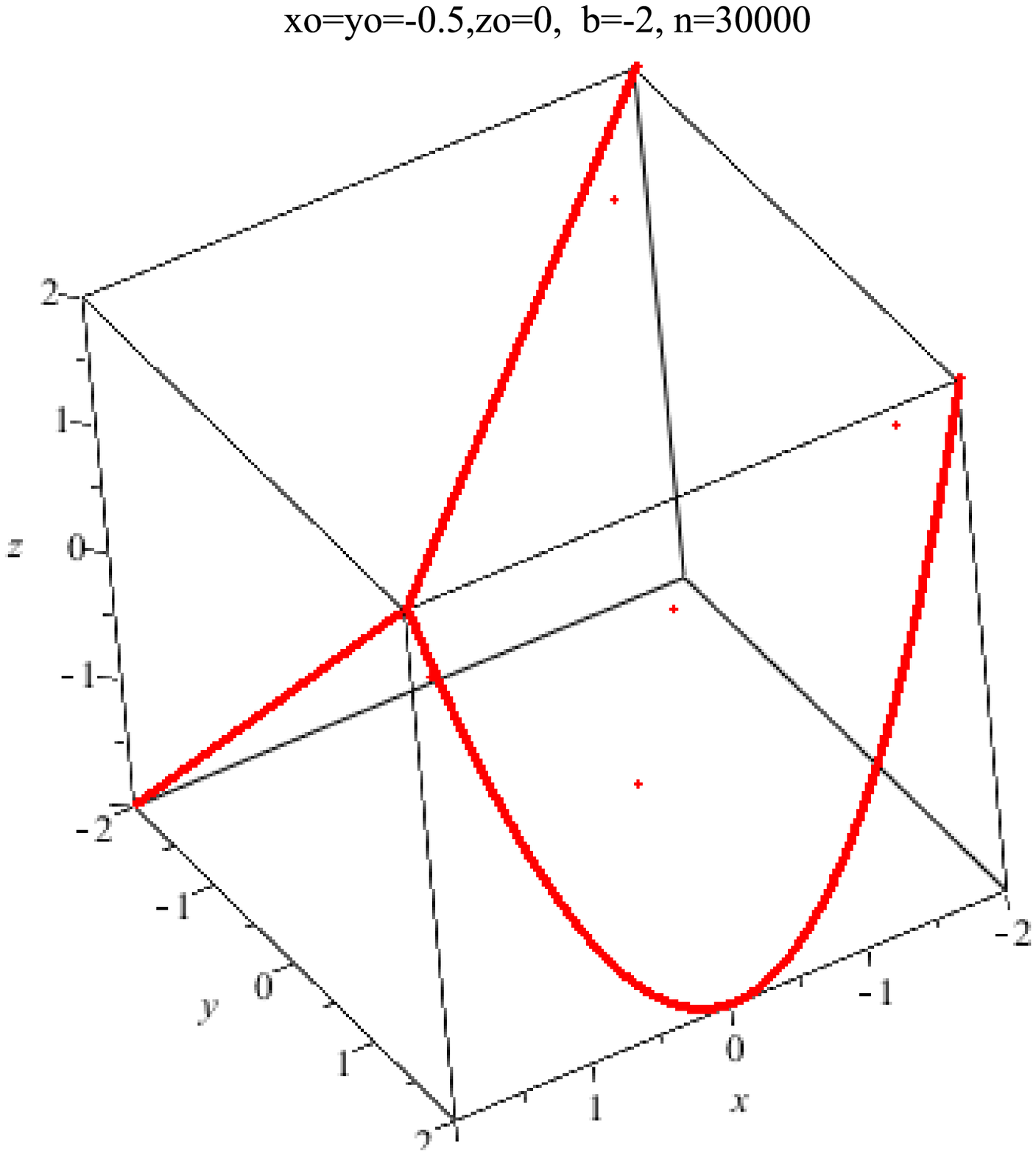}}
         \caption{Third order attractors for  different values of $X_{0}$ at $b=-2$.}\label{At5}
\end{center}
\end{figure}

The exponents of Lyapounov allowed us to conclude, for $b=-1.864$ (figure \ref{At4}) and $b=-2$ (figure \ref{At5}), that the attractors we obtained are chaotic attractors. Indeed, for $b = -1.864$, the Lyapounov exponents, calculated numerically, have the following values: $l_{1} = 0.1535$, $l_{2} = 0.1532$ and $l_{3} = 0.1532$, who are all positive. Similarly for $b=-2$, Lyapounov's exponents are: $l_{1} = 0.231048$, $l_{2} = 0.231048$ and $l_{3} = 0.231046$,  are all positive, which makes us allows to consider it chaotic.\\

\begin{figure}[!ht]
\begin{center}
    \subfigure[]{\label{At8-a}\includegraphics[width=3.6cm,height=3.5cm]{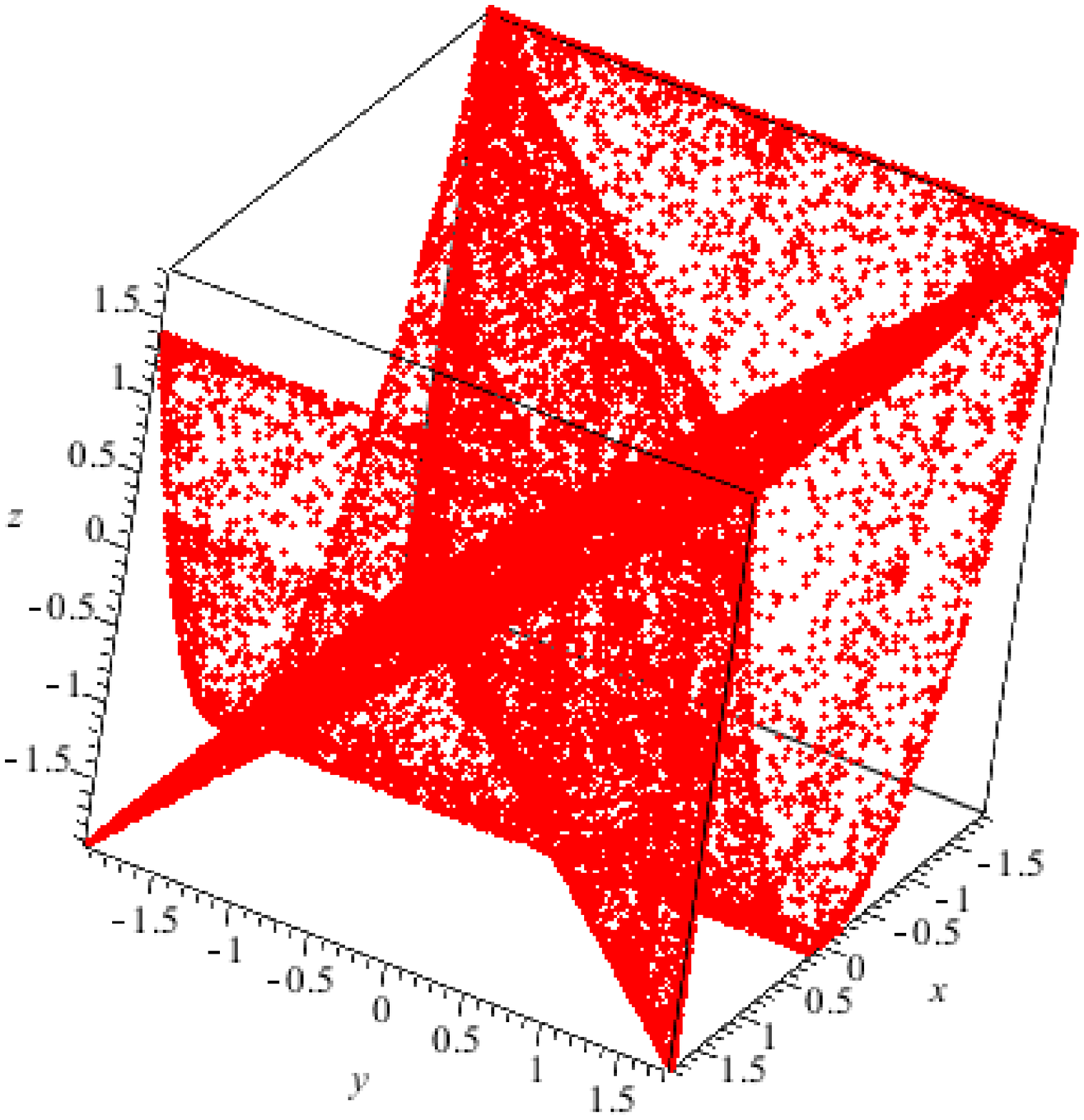}}
    \subfigure[]{\label{At7-c}\includegraphics[width=3.6cm,height=3.5cm]{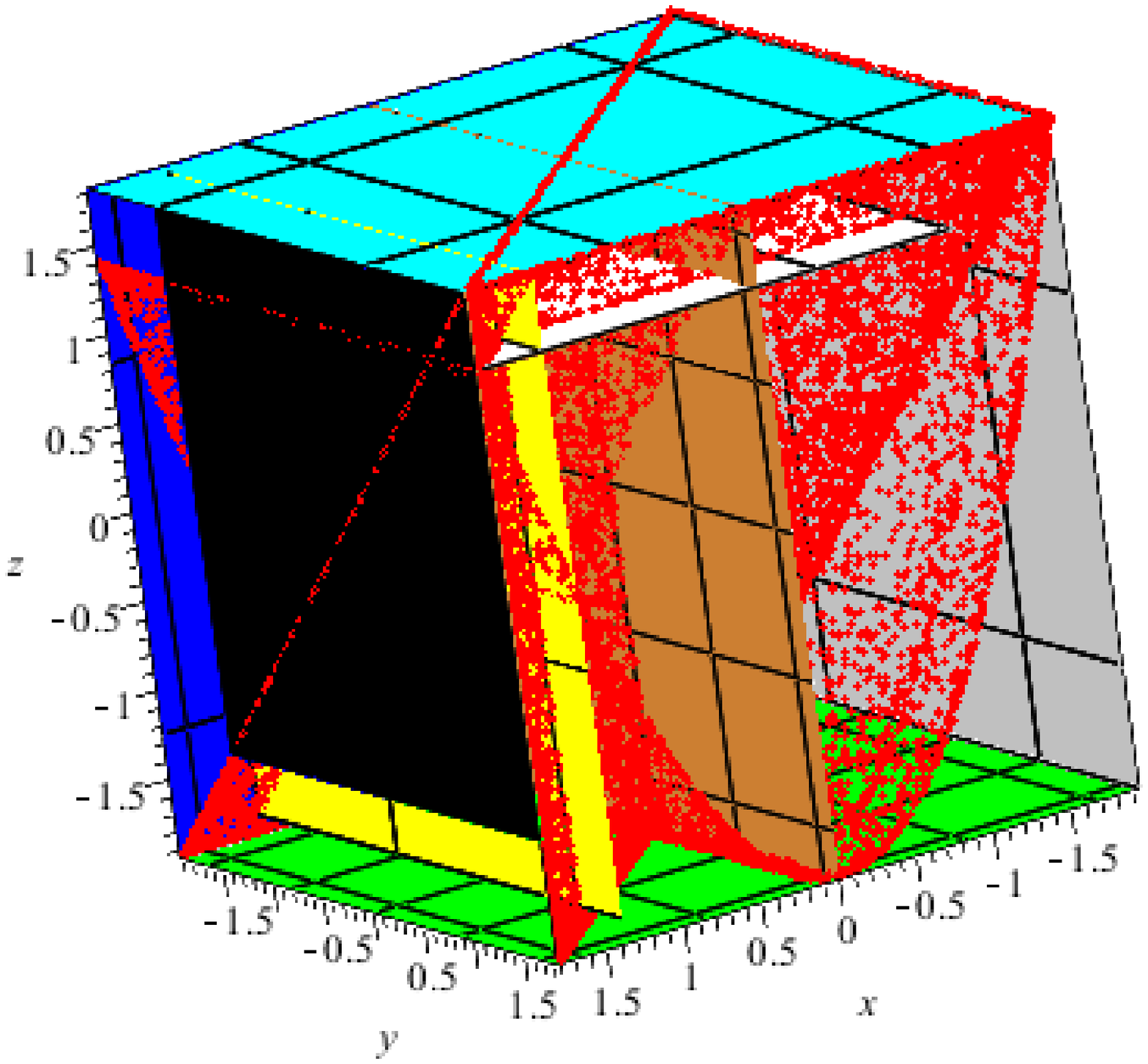}}
    \subfigure[]{\label{At7-d}\includegraphics[width=3.6cm,height=3.5cm]{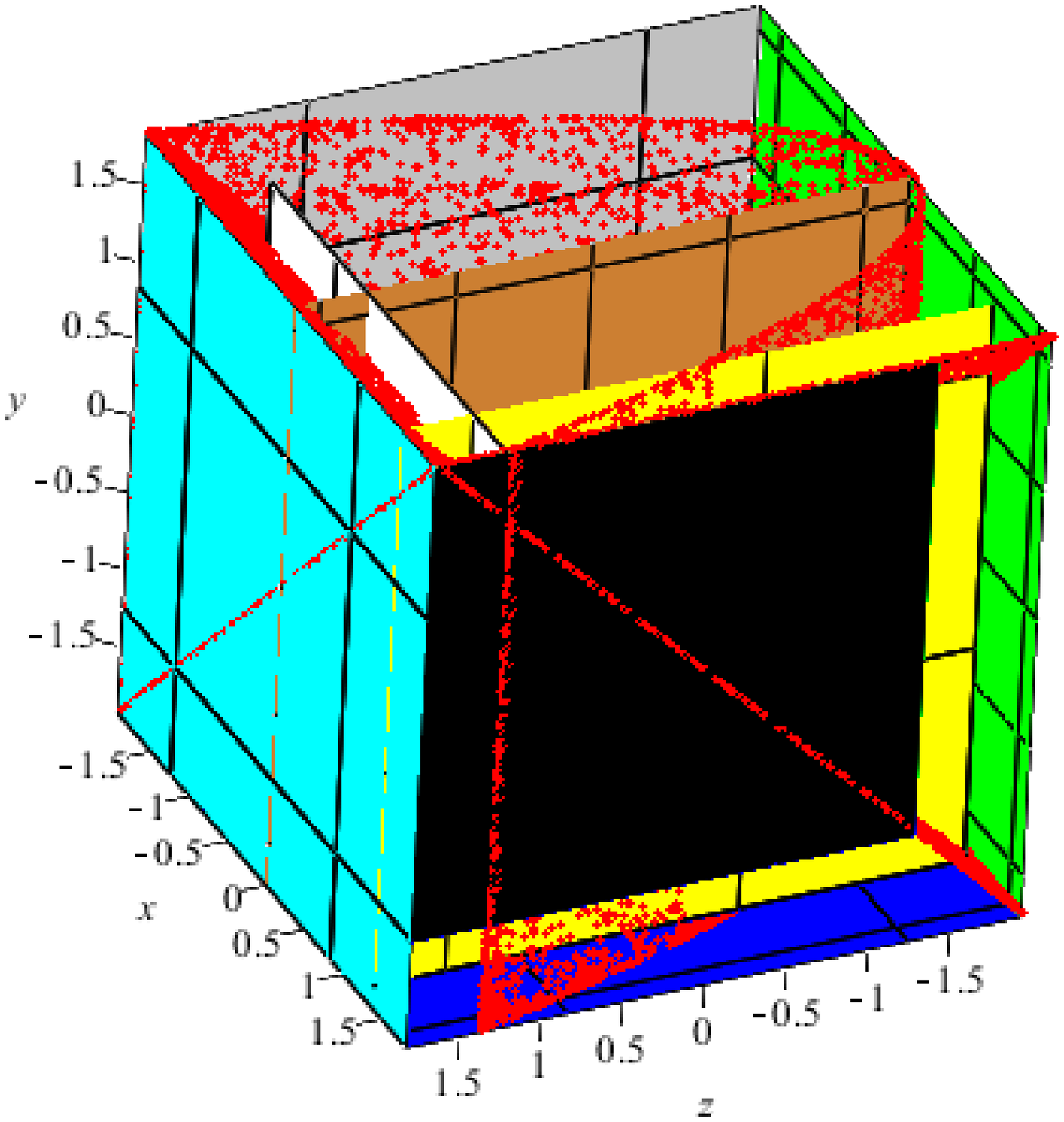}}
        \caption{Chaotic attractor bounded by critical planes at $b=-1.864$.}\label{At7}
\end{center}
\end{figure}
\begin{figure}[!ht]
\begin{center}
    \subfigure[]{\label{At8-a}\includegraphics[width=3.6cm,height=3.4cm]{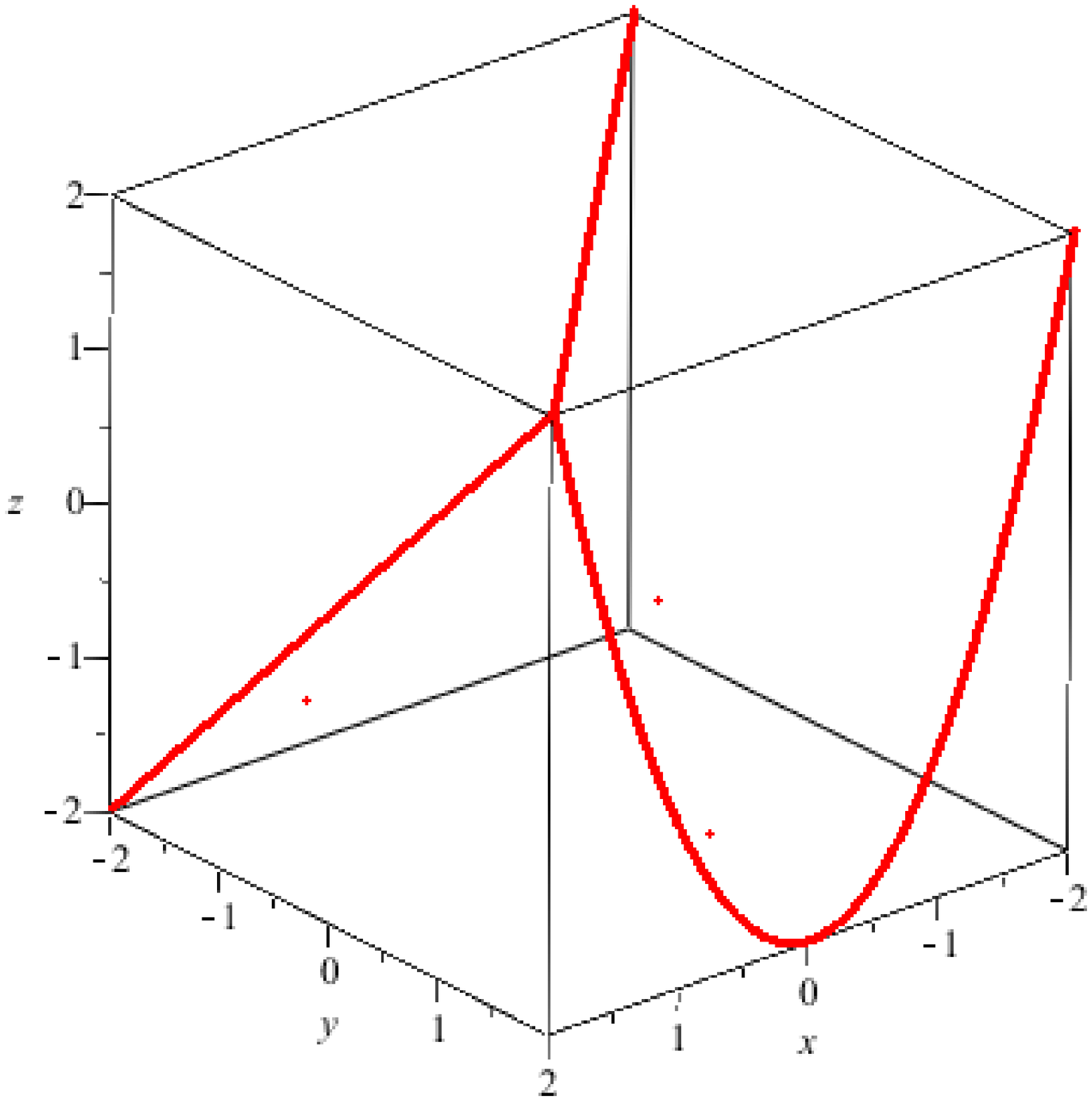}}
    \subfigure[]{\label{At8-b}\includegraphics[width=3.6cm,height=3.5cm]{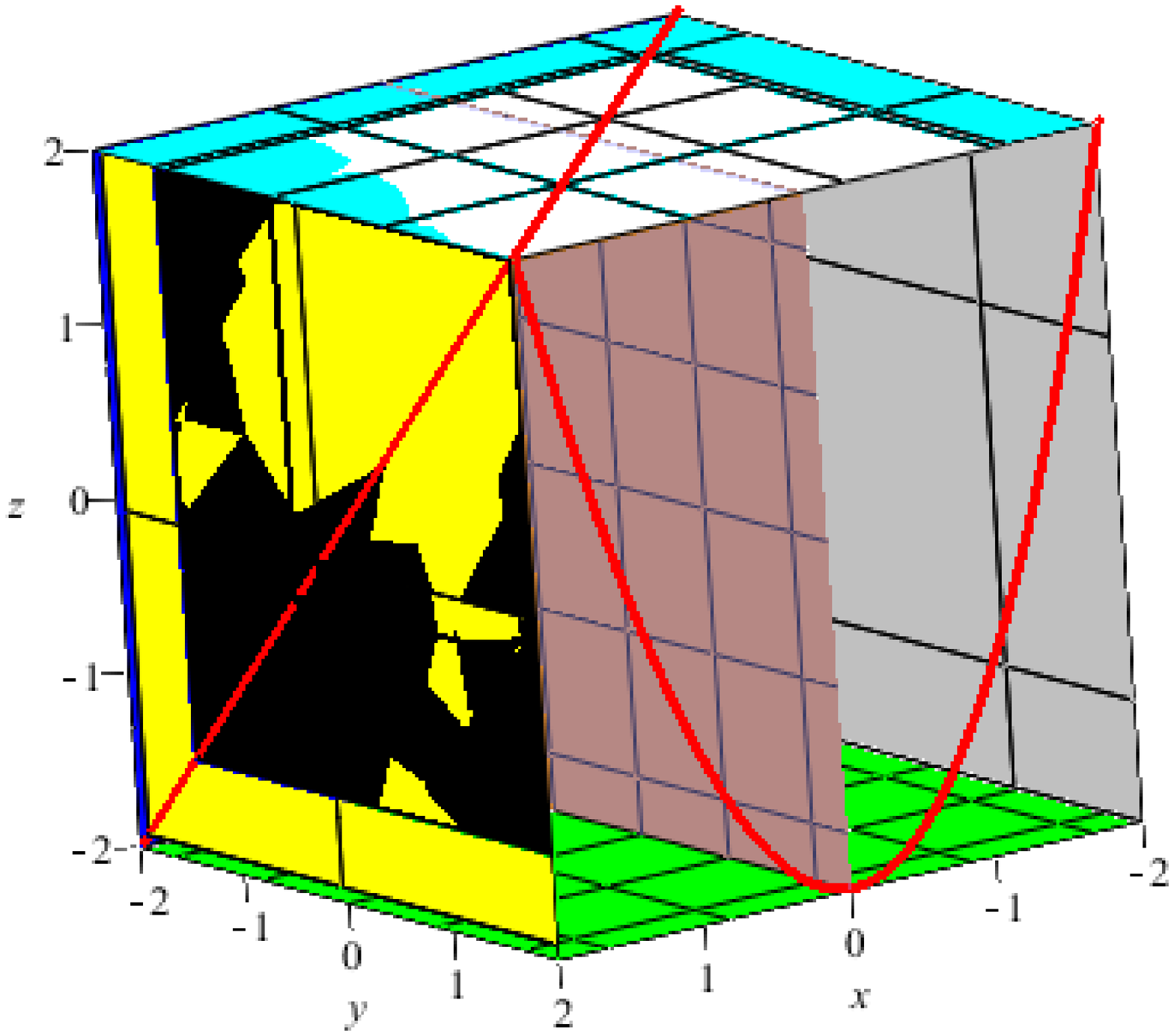}}
    \subfigure[]{\label{At8-c}\includegraphics[width=3.6cm,height=3.5cm]{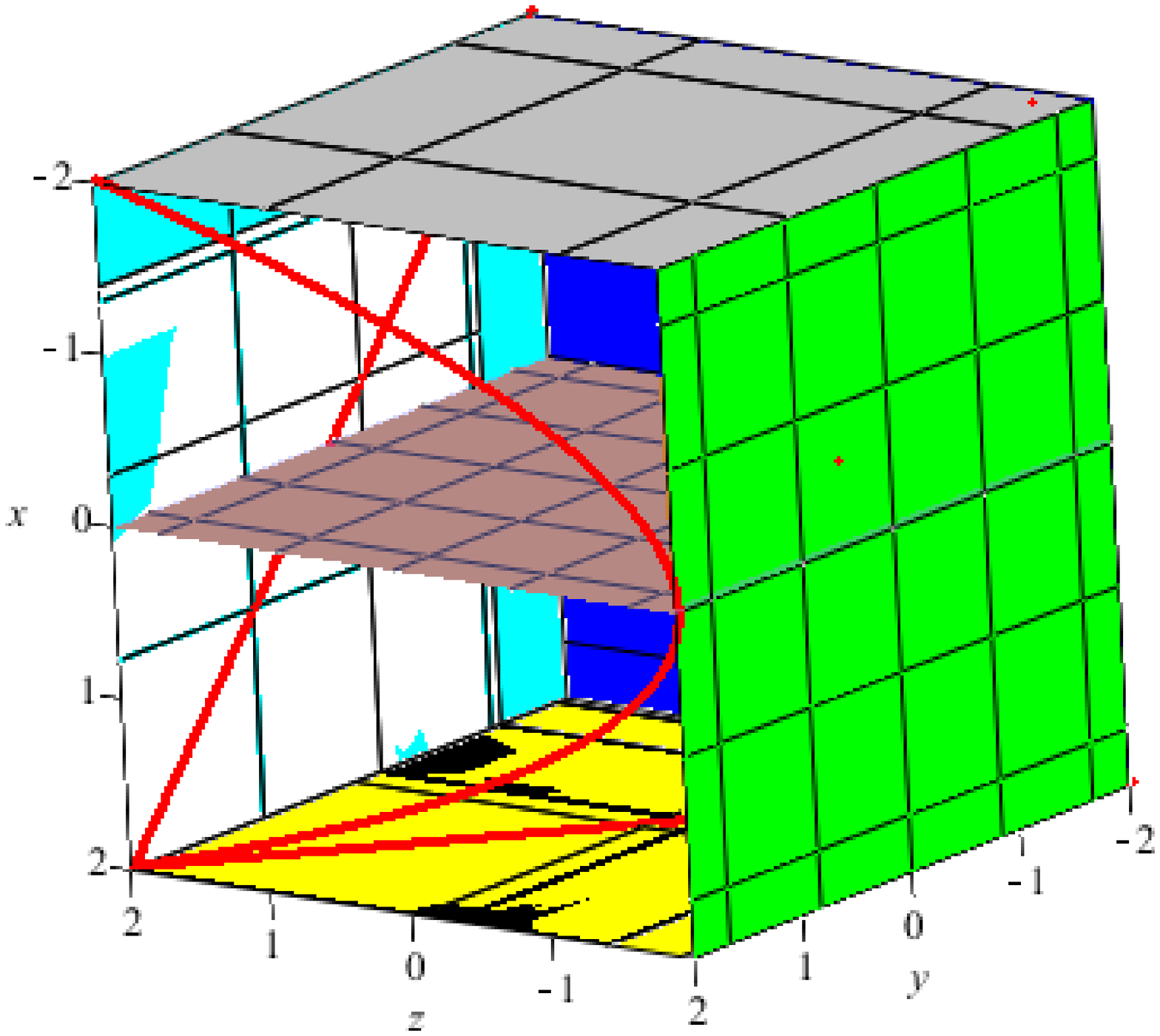}}
         \caption{Chaotic attractor bounded by critical planes at $b=-2$.}
  \label{At8}
\end{center}
\end{figure}
The figures \ref{At7} and \ref{At8} represent two orbits of $ 30 000 $ points, as well as some critical planes which delimit the chaotic attractors defined for $ b = -1.864 $ and $ b = -2 $ respectively, taking $ X_ {0} = (0, -0.5, 0.5) $ as the initial value.
Critical plans are represented in different colors: $PC_{-1}$ in brown, $PC$ in green, $PC_{1}$ in blue, $PC_{2}$ in gray, $PC_{3}$ in cyan, $PC_{5}$ in black, $PC_{6}$ in white and  $PC_{8}$ in yellow. \\
The chaotic attractors defined for these conditions correspond to unstable chaos. Chaos is said to be unstable when there is existence of a strange transient due to the presence of an infinity of unstable periodic solutions; we then speak of a chaotic repeller, such a set can be associated with the existence of an attractor at infinity (divergence for the initial conditions chosen) or the existence of a fuzzy boundary between the basins of two attractors. For $ b = -1.864 $, all fixed points are unstable nodes or unstable focus nodes. \\
\begin{figure}[!ht]
\begin{center}
\subfigure[$b=-1.864$]{\label{basAt1-a}\includegraphics[width=3.6cm,height=3.6cm]{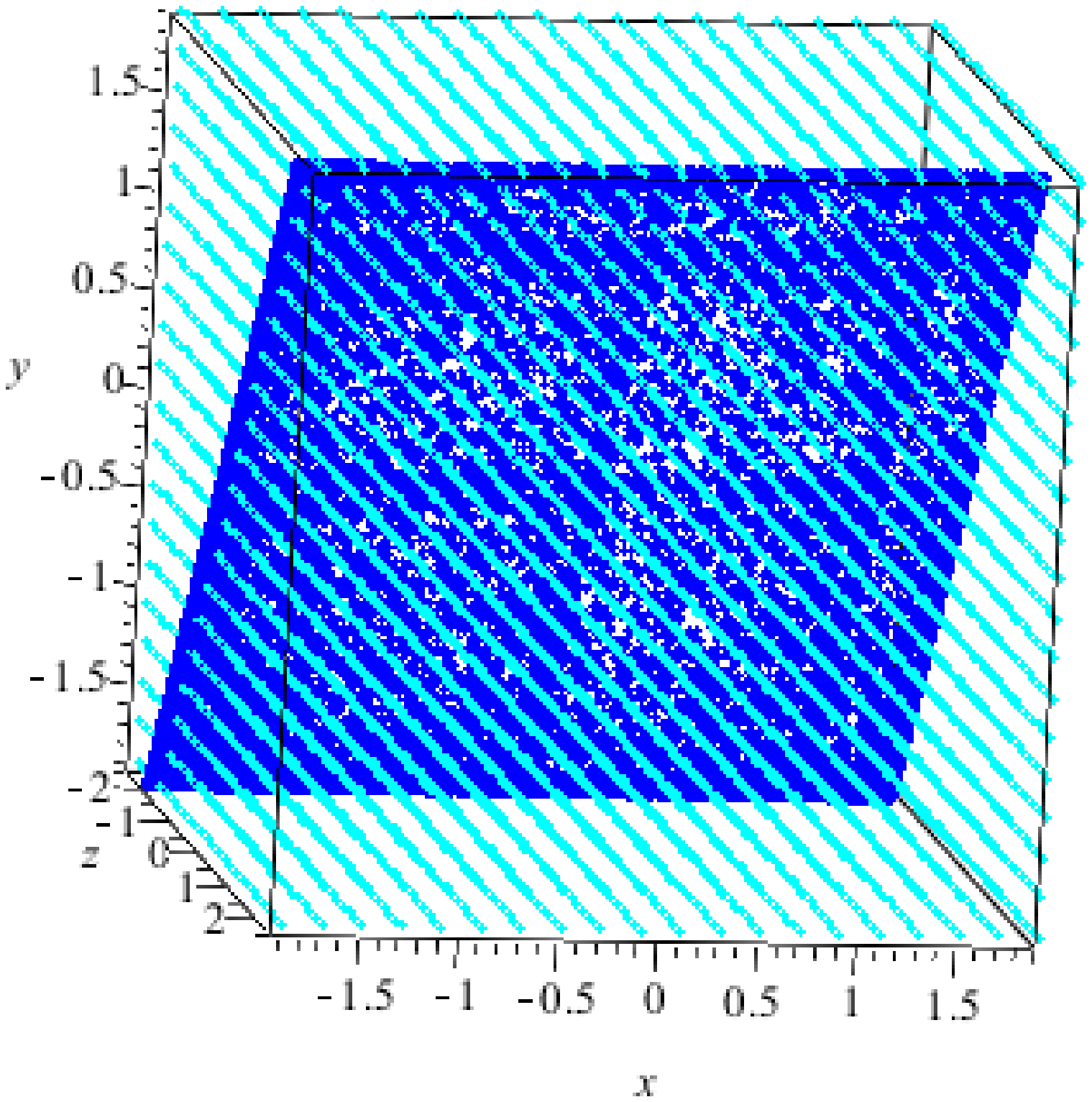}}
    \subfigure[$b=-2$]{\label{basAt1-b}\includegraphics[width=3.6cm,height=3.6cm]{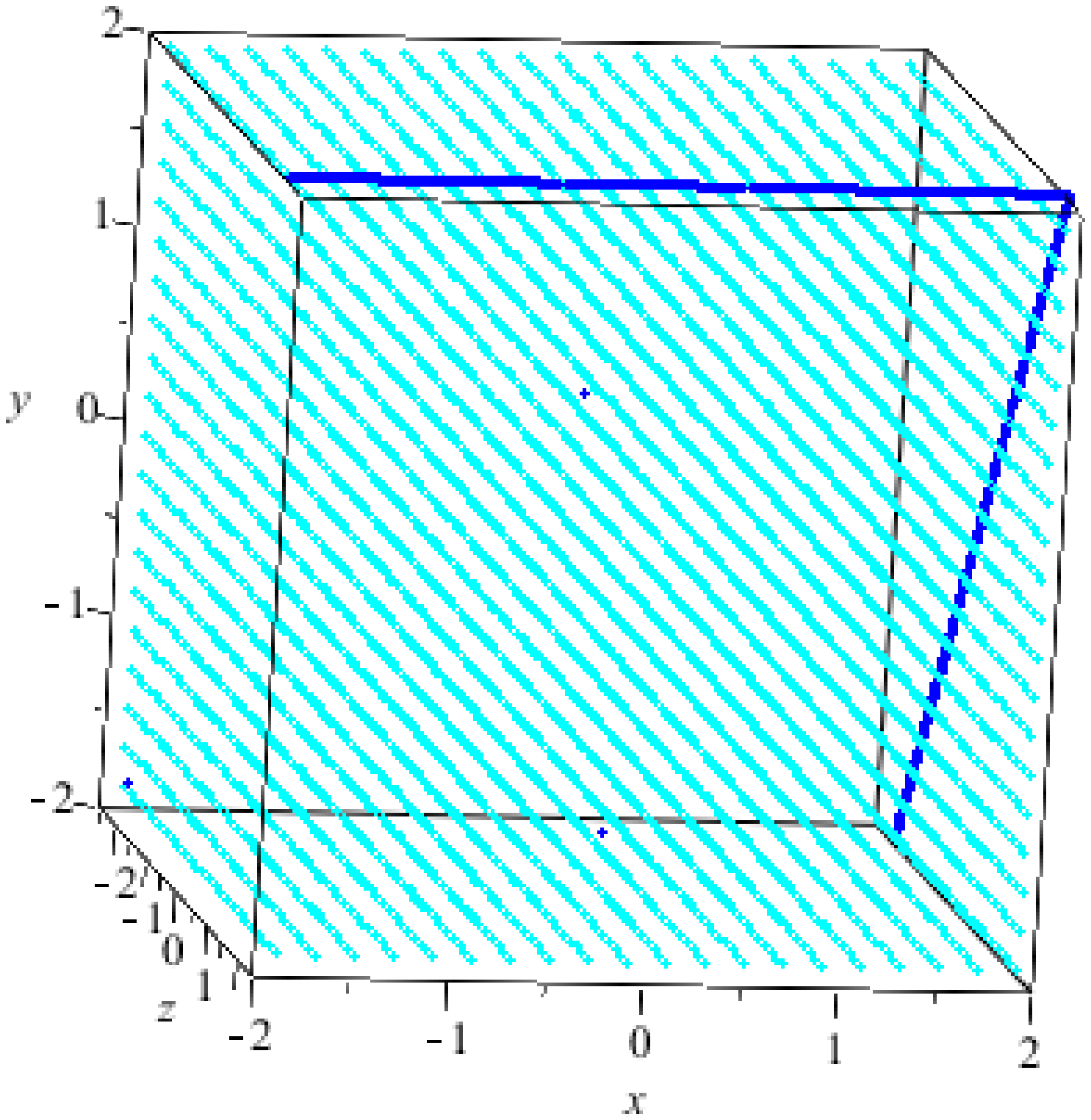}}
  \caption{Attraction Basins  with $X_{0}(0, -0.5, 0.5)$.}\label{basAt1}
  \end{center}
\end{figure}
 In the same way we consider the case of stable chaos, for which we have the existence of a neighborhood ensuring the convergence of the trajectories of the iterates. In other words, it is necessary to determine the initial values of the attractor basin allowing the convergence of the trajectories towards the attractor. The graphical representation of the attraction ponds of two dense orbits in the chaotic attractors defined for the two values of  $b$ ($b=-1.864$ and $b=-2$) with $X_{0}=(0, -0.5, 0.5)$, is given in the figure \ref{basAt1}.

\begin{figure}[!ht]
\begin{center}
    \subfigure[$b=-1.864$]{\label{BassinAtrachao8-a}\includegraphics[width=3.6cm,height=3.5cm]{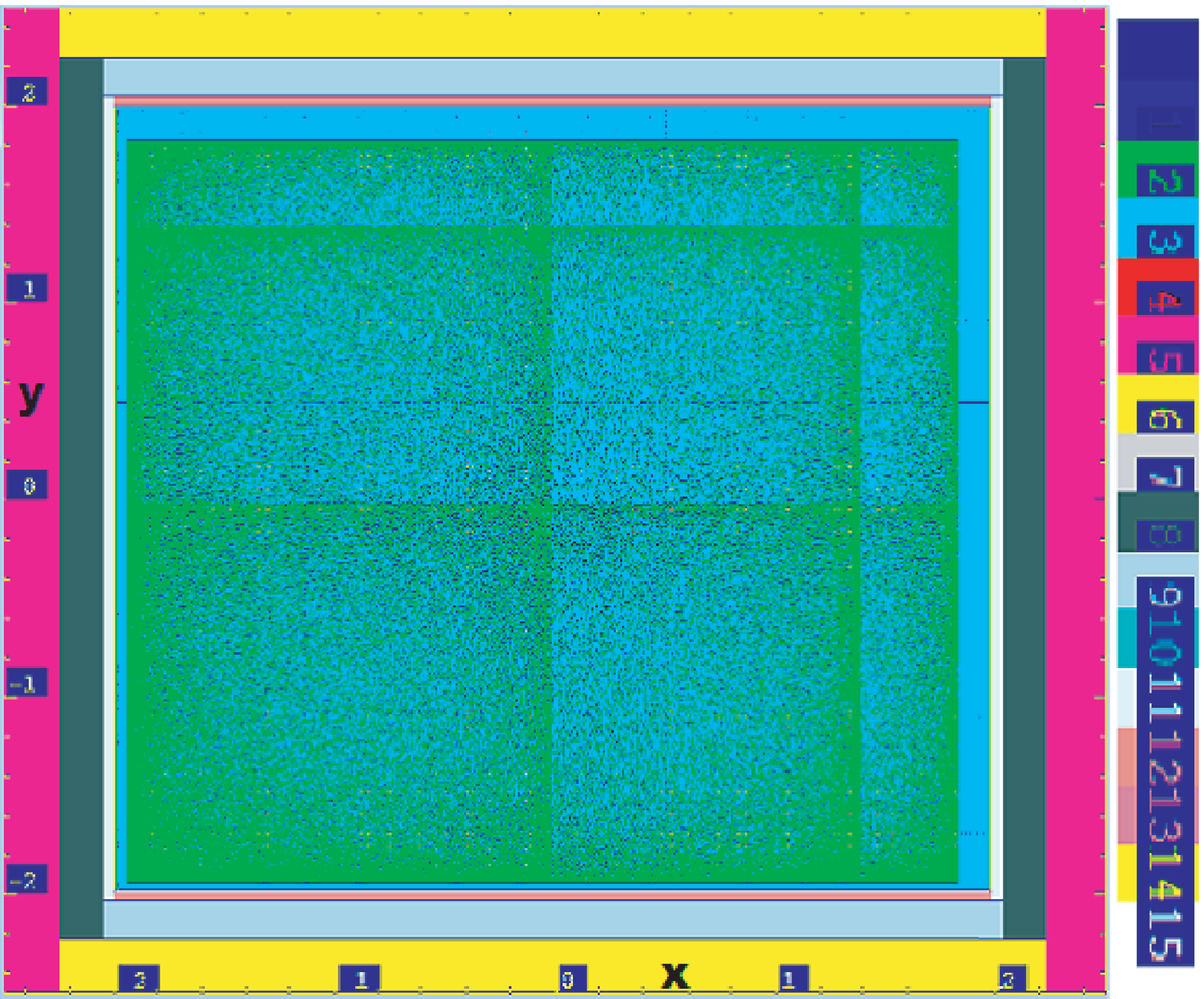}}
    \hspace{1cm}
    \subfigure[$b=-2$]{\label{BassinAtrachao8-b}\includegraphics[width=3.6cm,height=3.5cm]{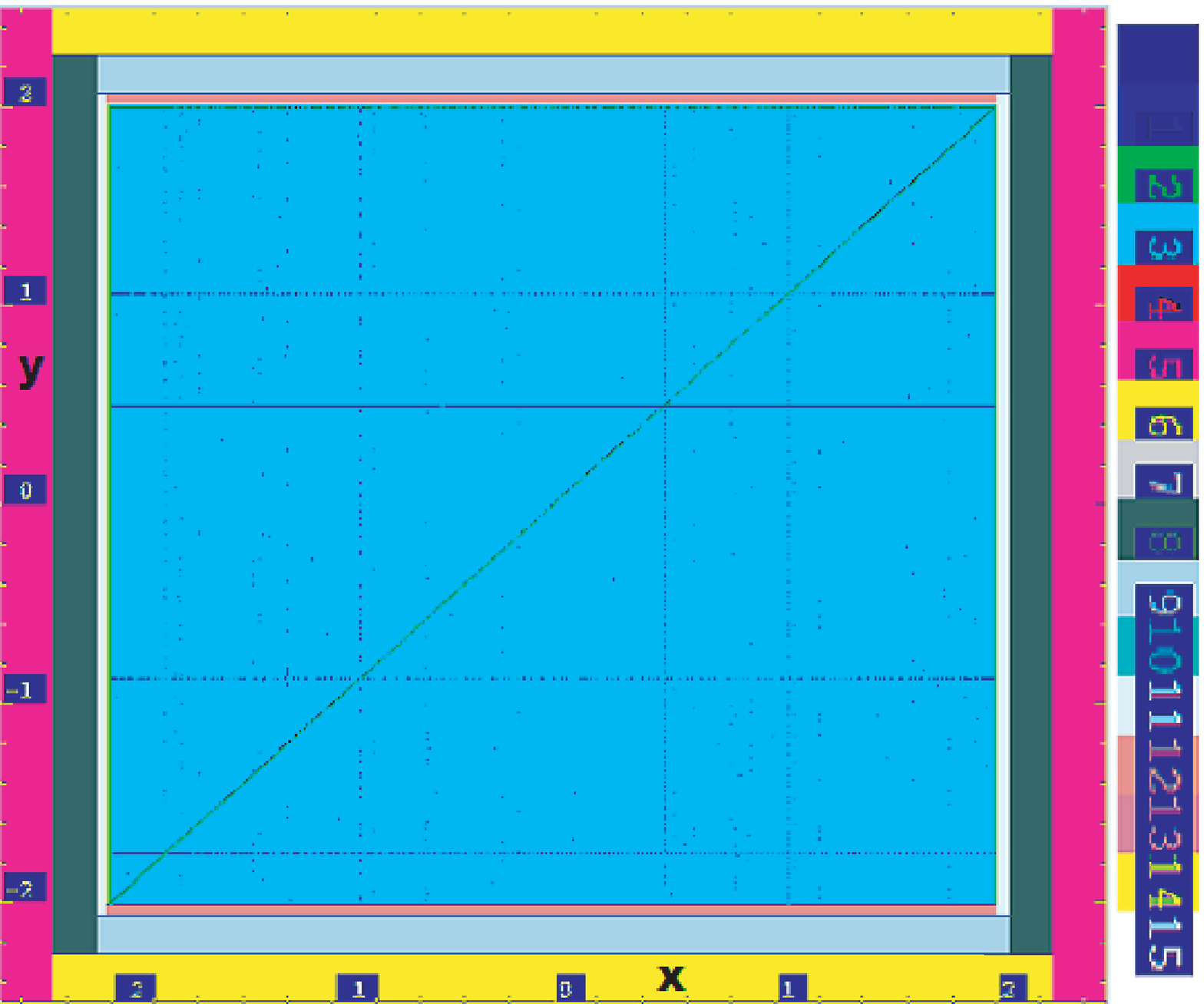}}\\
    \vspace{-0.3cm}
    \subfigure{\label{BassinAtrachao8-c}\includegraphics[width=6cm,height=0.6cm]{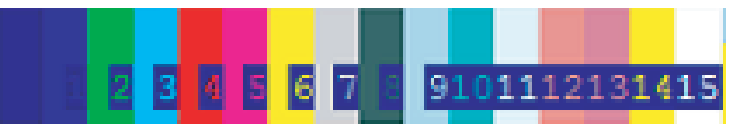}}
    \vspace{-0.3cm}
         \caption{Basins of chaotic attractors of $T$, for two values of $b$ ($X_{0}(0, -0.5, 0.5)$).}
  \label{BassinAtrachao8}
\end{center}
\end{figure}

 We note that the two orbits touch the boundaries of their basins of attraction. For the same values of $ b $, we have
  the existence of several chaotic attractors as shown in the figure \ref{BassinAtrachao8}, whose basins of attraction of coexisting attractors are represented with different colors.

\end{document}